\documentclass[11pt]{article}

\usepackage[latin1]{inputenc}

\usepackage{amsfonts,amsmath,amssymb,amsthm,graphicx,epsfig,float}
\usepackage[T1]{fontenc}

\usepackage[english,francais]{babel}

{
  \newtheorem{theoreme}{Th\'eor\`eme}
  \newtheorem*{theoreme*}{Th\'eor\`eme}
  \newtheorem{lemme}[theoreme]{Lemme}
  \newtheorem{corollaire}[theoreme]{Corollaire}

  \newtheorem{proposition}[theoreme]{Proposition}
\newtheorem*{corollaire*}{Corollaire}
\newtheorem*{proposition*}{Proposition}
\theoremstyle{remark}
  \newtheorem*{remarque*}{Remarque}
}

\newcounter{ex}

\newenvironment{rem*}{
  \noindent\textbf{Remarque. }}{}

\newcommand{\Cc}{\mathbb{C}}
\newcommand{\Nn}{\mathbb{N}}
\newcommand{\Rr}{\mathbb{R}}

\newcommand{\Zz}{\mathbb{Z}}
\newcommand{\Pp}{\mathbb{P}}
\newcommand{\Acal}{\mathcal{A}}

\newcommand{\Hcal}{\mathcal{H}}
\newcommand{\Ical}{\mathcal{I}}

\newcommand{\Scal}{\mathcal{S}}
\newcommand{\Wcal}{\mathcal{W}}

\title{{\bf Sur les exposants de Lyapounov des applications méromorphes}}
\author{Henry de Thélin}
\date{}

\begin{document}
\maketitle


\def\figurename{{Fig.}}%
\def\proofname{Preuve}
\def\contentsname{Sommaire}%

\begin{abstract}

Soit $f$ une application méromorphe dominante d'une variété
Kählérienne compacte. Nous donnons une inégalité pour les exposants de
Lyapounov d'une classe de mesures ergodiques de $f$ en utilisant
l'entropie métrique et les degrés dynamiques de $f$. Nous en déduisons
l'hyperbolicité de certaines mesures.

\end{abstract}

\selectlanguage{english}
\begin{center}
{\bf{ }}
\end{center}

\begin{abstract}

Let $f$ be a dominating meromorphic self-map of a compact Kähler
manifold. We give an inequality for the Lyapounov exponents of some ergodic
measures of $f$ using the metric entropy and the
dynamical degrees of $f$. We deduce the hyperbolicity of some measures.

\end{abstract}

\selectlanguage{francais}

Mots-clefs: applications méromorphes, exposants de Lyapounov,
entropie.

Classification: 37Fxx, 32H50, 58F15.

\section*{{\bf Introduction}}
\par

Soit $(X, \omega)$ une variété Kählérienne compacte de dimension $k$ et
$f : X \mapsto X$ une application méromorphe dominante.

Nous désignerons par $C_f$ l'ensemble critique de $f$ et par $I_f$ son
ensemble d'indétermination.

L'objet de cet article est de donner des formules générales pour les
exposants de Lyapounov des mesures invariantes $\mu$ qui intègrent la fonction
$\log d(x ,\Acal)$ où $d$ est la distance dans $X$ et $\Acal=C_f \cup
I_f$. Remarquons que lorsqu'une mesure $\mu$ intègre la fonction précédente,
elle ne charge pas l'ensemble $\Acal$: on peut donc définir $f_{*}
\mu$ et parler de mesure invariante. Par ailleurs, l'hypothèse
d'intégrabilité de $\log d(x ,\Acal)$ est vérifiée dès que $\mu$ intègre les fonctions quasi-psh. Les
formules dépendront d'une part de l'entropie métrique $h(\mu)$ de
$\mu$ et d'autre part des degrés dynamiques $d_q$ de $f$ (voir le
paragraphe \ref{rappel} pour leur définition).

Dans ce contexte nous avons le

\begin{theoreme}{\label{formule}}

Soient $\mu$ une mesure invariante, ergodique telle que $\log d(x
,\Acal) \in L^1(\mu)$ et $\chi_1 \geq \dots \geq \chi_k$ les exposants
de Lyapounov de $\mu$ (ils sont bien définis).

Fixons $1 \leq s \leq k$. On définit $l=l(s)$ et $l'=l'(s)$ par
les formules suivantes:

$$\chi_1 \geq \dots \geq \chi_{s-l-1} > \chi_{s-l}= \dots = \chi_s =
\dots = \chi_{s+l'} > \chi_{s+l'+1} \geq \dots \geq \chi_k,$$
où $s-l$ est égal à $1$ si $\chi_1= \dots = \chi_s$ et $s+l'$ est égal
à $k$ lorsque $\chi_s = \dots = \chi_k$.

Alors, on a les inégalités suivantes:

$$h(\mu) \leq \max_{0 \leq q \leq s-l-1} \log d_q +2 \chi_{s-l}^{+} +
\dots + 2 \chi_{k}^{+}$$

$$h(\mu) \leq \max_{s+l' \leq q \leq k} \log d_q -2 \chi_1^{-} - \dots -
2 \chi_{s+l'}^{-}$$

avec $\chi_i^{+}= \max ( \chi_i, 0)$ et $\chi_i^{-}= \min (\chi_i,0)$.

\end{theoreme}

Signalons la ressemblance entre ces formules et celles de J. Buzzi
pour les applications $C^{1 + \alpha}$ (voir \cite{Bu}).

Maintenant, à l'aide de notre théorème, on a:

\begin{corollaire}{\label{cor1}}

Supposons que les degrés dynamiques vérifient  $d_1 \leq
\dots \leq d_{s-1} < d_s > d_{s+1} \geq \dots \geq d_k$. Soit $\mu$ une mesure invariante, ergodique telle que $\log d(x
,\Acal) \in L^1(\mu)$ et $h(\mu) > \max ( \log d_{s-1}, \log
d_{s+1})$ (ou $h(\mu) > \log d_{k-1}$ si $s=k$). Alors
$$\chi_1 \geq \dots \geq \chi_s \geq \frac{1}{2}(h(\mu) - \log d_{s-1}) > 0$$
et
$$0 >  \frac{1}{2} ( \log d_{s+1} - h(\mu)) \geq \chi_{s+1} \geq \dots \geq \chi_k.$$

En particulier la mesure $\mu$ est hyperbolique.

\end{corollaire}

Lorsque dans le corollaire précédent la mesure $\mu$
est d'entropie $\log d_{s}$ (i.e. est d'entropie maximale par
\cite{DS1} et \cite{DS2}), on obtient le:
 
\begin{corollaire}{\label{cor2}}

Soit $\mu$ une mesure invariante, ergodique telle que $\log d(x
,\Acal) \in L^1(\mu)$. Alors, si $h(\mu) = \log d_s$ avec $d_1 \leq
\dots \leq d_{s-1} < d_s > d_{s+1} \geq \dots \geq d_k$ on a:
$$\chi_1 \geq \dots \geq \chi_s \geq \frac{1}{2} \log
\frac{d_s}{d_{s-1}} > 0$$
et
$$0 > \frac{1}{2} \log
\frac{d_{s+1}}{d_{s}} \geq \chi_{s+1} \geq \dots \geq \chi_k.$$
En particulier la mesure $\mu$ est hyperbolique.

\end{corollaire} 

Notons que les inégalités obtenues dans ce corollaire sont celles qui étaient
conjecturées (voir \cite{Gu2} Conjecture 3.2).

Les hypothèses de ce corollaire sont vérifiées dans de nombreuses
situations. En voici certaines.

Tout d'abord pour les endomorphismes holomorphes de
$\Pp^k$ avec la mesure de Green $\mu$ (voir \cite{FS} et
\cite{FS1} pour sa définition). En effet, la mesure $\mu$ est
mélangeante et intègre $\log d(x,\Acal)$ grâce
à l'inégalité de Chern-Levine-Nirenberg. Pour ces endomorphismes,
l'entropie métrique de $\mu$ vaut $k \log(d)$ et les $d_q$ valent
$d^q$. En appliquant notre inégalité avec $s=k$, on a alors la
minoration du plus petit exposant de $\mu$ par $\frac{\log(d)}{2}$ et
on retrouve ainsi un résultat de J.-Y. Briend et J. Duval (voir
\cite{BD1}).

De la même façon, lorsque $f$ est une application méromorphe sur une
variété projective, avec son degré topologique strictement plus grand
que les autres (i.e. $s=k$ dans le corollaire précédent), V. Guedj a construit
une mesure $\mu$ mélangeante, d'entropie $\log d_k$ et qui intègre les
fonctions quasi-psh (voir \cite{Gu} et \cite{DS1}). En utilisant notre formule, on retrouve alors la
minoration du plus petit exposant de Lyapounov de $\mu$ par
$\frac{1}{2} \log(d_k / d_{k-1})$ qu'il avait démontrée.

Lorsque $f$ est un automorphisme holomorphe d'une variété Kählérienne qui possède un degré dynamique strictement plus grand que les
autres, T.-C. Dinh et N. Sibony ont construit une mesure de Green
$\mu$ mélangeante qui intègre les fonctions quasi-psh et d'entropie
$\max \log d_q$ (voir \cite{DS3}). Notre corollaire s'applique donc et
on obtient ainsi un nouveau résultat pour ces automorphismes.

De la même façon, lorsque $f$ est une application birationnelle
régulière de $\Pp^k$ (voir \cite{DS4}), T.-C. Dinh et N. Sibony ont construit une mesure de Green
$\mu$ mélangeante et qui intègre la fonction $\log d(x,\Acal)$ (ici
$I_f=I^{+}$ et $C_f= f^{-1}(I^{-})$). On peut donc lui appliquer le corollaire.

Donnons une autre conséquence de notre théorème. Lorsque l'on applique la première formule du théorème avec $s=1$ on en
déduit une formule de Ruelle (\cite{Ru}) pour les applications méromorphes:

\begin{corollaire}

Soit $\mu$ une mesure invariante, ergodique telle que $\log d(x
,\Acal) \in L^1(\mu)$. Alors:

$$h(\mu) \leq 2 \chi_1^{+} + \dots + 2 \chi_k^{+}.$$

\end{corollaire}

Voici le plan de ce texte. Dans un premier paragraphe nous ferons des
rappels sur les applications méromorphes et dans le second nous
démontrerons les corollaires \ref{cor1} et \ref{cor2}. Dans le troisième nous parlerons de
théorie de Pesin pour les applications méromorphes et dans le
quatrième, nous ferons des rappels sur la transformée de
graphe. Enfin, le cinquième paragraphe sera consacré à la
démonstration de la première inégalité du théorème et le sixième à celle
de la deuxième inégalité. Dans le dernier paragraphe nous donnerons un
analogue de notre théorème pour les difféomorphismes de classe $C^{1 +
  \alpha}$ dans les variétés Riemanniennes compactes.

\section{{\bf Rappels sur les applications méromorphes}}{\label{rappel}}

Commençons par rappeler la définition des degrés dynamiques de $f$
(voir \cite{RS} et \cite{DS2}).

On définit la forme $f^{*} (\omega^q)$ comme l'extension triviale de $(f_{|X
  \setminus I_f})^{*} \omega^q$.

On pose $\delta_q(f):=
  \int_X f^{*}( \omega^q) \wedge \omega^{k-q}$ pour $q=0, \dots, k$. Le
  degré dynamique d'ordre $q$ est alors:
$$d_q:= \lim_{n \rightarrow \infty} ( \delta_q(f^n))^{1/n}.$$

Notons que $d_0=1$ et que la limite ci-dessus existe grâce à
\cite{DS2}.

Dans \cite{DS1} et \cite{DS2}, T.-C. Dinh et N. Sibony ont défini, via
un procédé de régularisation de courant, le pull-back par $f$ des
courants positifs fermés de bidimension quelconque. En voici le
procédé. Tout d'abord $f$ est une submersion sur un ouvert de Zariski $\Omega_{1,f}$ de $X$. Si $S$ est un courant positif fermé de bidegré
$(l,l)$ sur $X$, alors, on peut définir $f^{*}S$ sur $\Omega_{1,f}$. C'est
un courant positif fermé dont la masse $\int_{\Omega_{1,f}} f^{*}S \wedge
\omega^{k-l}$ est majorée par $c_X
\delta_l(f) \| S \|$ (voir le lemme 4 de \cite{DS1} et le
corollaire 1.3 de \cite{DS2}, ici $c_X$ est une
constante qui ne dépend que de $X$). Ce courant admet donc un prolongement
trivial $\widetilde{f^{*}S}$ à $X$ tout entier d'après un théorème de
H. Skoda (\cite{S}). De plus $\widetilde{f^{*}S}$ est un courant positif
fermé de masse majorée par $c_X \delta_l(f) \| S \|$. Ce courant sera
appelé le pull-back de $S$ par $f$.

Nous allons maintenant donner deux lemmes qui serons utilisés dans la
démonstration du théorème. Le premier est quasiment le même que le
Lemme 5 de \cite{DS1}. On notera $\Omega_f= X \setminus
\cup_{n \in \Zz} f^n(I_f)$.

\begin{lemme}{\label{lemme1}}

Soit $\epsilon > 0$. Il existe une constante $c_{\epsilon} > 0$ telle
que pour tout $q=0, \dots , k$ on ait:

$$\int_{\Omega_f } \omega^{k-q} \wedge (f^{n_1})^{*} \omega \wedge
\dots \wedge (f^{n_q})^{*} \omega \leq c_{\epsilon} ( \max_{0 \leq j
  \leq q }d_j + \epsilon)^{n_1}$$
pour tous les entiers naturels $n_1 \geq \dots \geq n_q \geq 0$.

\end{lemme}

\begin{proof}

La démonstration est la même que dans \cite{DS1}. Nous la
donnons par confort pour le lecteur.

Soit $c > 0$ une constante telle que $\delta_j(f^n) \leq c (d_{j} +
\epsilon)^n$ pour tout $n \geq 0$ et tout $j=0, \dots , k$.

Soit $\Omega_{n,f} = X \setminus \cup_{0 \leq i \leq n-1} f^{-i}(C)$
où $C= X \setminus \Omega_{1,f}$. On va montrer par
récurrence sur $q$, avec $0 \leq q \leq k$ que pour tous $n_1 \geq
\dots \geq n_q \geq 0$, on a $\| T_q \|= \int_{\Omega_{n_1 ,f}} T_q
  \wedge \omega^{k-q}  \leq c^q c_{X}^q (\max_{0 \leq j \leq q} d_j +
\epsilon)^{n_1}$ où
$$T_q=(f^{n_1})^{*} \omega \wedge \dots \wedge (f^{n_q})^{*} \omega
,$$
et $T_0=1$.

C'est vrai pour $q=0$. Supposons la propriété vraie au rang $q-1$. Cela
implique que $\| T_{q-1}' \| \leq c^{q-1} c_{X}^{q-1} (\max_{0 \leq j
  \leq q-1} d_{j} +
\epsilon)^{n_1-n_q}$ avec
$$T_{q-1}'=(f^{n_1-n_q})^{*} \omega \wedge \dots \wedge (f^{n_{q-1} -
  n_q})^{*} \omega .$$
Le courant $T_{q-1}'$ est donc de masse finie sur
  $\Omega_{n_1-n_q,f}$. Il admet donc une extension triviale
  $\widetilde{T_{q-1}'}$ dans $X$ qui est un courant positif fermé de
  masse majorée par $c^{q-1} c_{X}^{q-1} (\max_{0 \leq j \leq q-1} d_{j} +
\epsilon)^{n_1-n_q}$. En utilisant la propriété du pull-back de
  T.-C. Dinh et N. Sibony énoncée avant le lemme avec
  $S=\widetilde{T_{q-1}'} \wedge \omega$, on obtient:
$$ \| T_q \|= \| (f^{n_q})^{*}(T_{q-1}' \wedge \omega) \| \leq c_X
    \delta_{q}(f^{n_q}) \| T_{q-1}' \|   \leq  c^q c_{X}^q (\max_{1
    \leq j \leq q} d_j + \epsilon)^{n_1}.$$
Cela démontre bien le lemme.

\end{proof}

On utilisera aussi le lemme suivant:

\begin{lemme}{\label{lemme2}}

Soit $\epsilon > 0$. Il existe une constante $c_{\epsilon} > 0$ telle
que pour tout $q=0 , \dots , k$ on ait:

$$\int_{\Omega_f } (f^n)^{*} \omega^{q} \wedge (f^{n_1})^{*} \omega \wedge
\dots \wedge (f^{n_{k-q}})^{*} \omega \leq c_{\epsilon} ( \max_{q \leq j
  \leq k }d_j + \epsilon)^{n}$$
pour tous les entiers naturels $n \geq n_1 \geq \dots \geq n_{k-q} \geq 0$.

\end{lemme}

\begin{proof}

Soit $c > 0$ une constante telle que $\delta_j(f^n) \leq c (d_{j} +
\epsilon)^n$ pour tout $n \geq 0$ et tout $j=0, \dots , k$. Fixons $q$
entre $0$ et $k$. On peut supposer $q >0$ sinon on a le résultat par
le lemme précédent. 

On va montrer par récurrence sur $j$, avec $0 \leq j \leq k-q$ que
pour tous $n \geq n_1 \geq
\dots \geq n_j \geq 0$, on a $\| T_j \|= \int_{\Omega_{n ,f}} T_j
  \wedge \omega^{k-j-q}  \leq c^{j+1} c_{X}^j (\max_{q \leq l \leq q+j} d_l +
\epsilon)^{n}$ où
$$T_j=(f^n)^{*} \omega^q \wedge (f^{n_1})^{*} \omega \wedge \dots \wedge (f^{n_j})^{*} \omega
,$$
et $T_0= (f^n)^{*} \omega^q$.

C'est vrai pour $j=0$ par définition du $q$-ème degré dynamique. Supposons la propriété vraie au rang $j-1$. Cela
implique que $\| T_{j-1}' \|  \leq c^{j} c_{X}^{j-1} (\max_{q \leq l \leq q+j-1} d_l +
\epsilon)^{n-n_j}$ où
$$T_{j-1}'= (f^{n-n_j})^{*} \omega^q \wedge (f^{n_1-n_j})^{*} \omega
\wedge \dots \wedge (f^{n_{j-1}-n_j})^{*} \omega.$$
Le courant $T_{j-1}'$ est donc de masse finie sur
  $\Omega_{n-n_j,f}$. Il admet donc une extension triviale
  $\widetilde{T_{j-1}'}$ dans $X$ qui est un courant positif fermé de
  masse majorée par $c^{j} c_{X}^{j-1} (\max_{q \leq l \leq q+j-1} d_l +
\epsilon)^{n-n_j}$. En utilisant encore la propriété du pull-back de
T.-C. Dinh et N. Sibony avec $S=\widetilde{T_{j-1}'} \wedge \omega$, on obtient:

$$\| T_j \| = \| (f^{n_j})^{*}(T_{j-1}' \wedge \omega) \| \leq c_{X}
  \delta_{q+j}(f^{n_j}) \|T_{j-1}'\| \leq c^{j+1} c_X^j (\max_{q \leq l \leq q+j} d_l +
\epsilon)^{n}.$$

Cela démontre la récurrence et quand on prend $j= k-q$, on obtient le lemme.

\end{proof}

\section{{\bf Démonstration des corollaires \ref{cor1} et \ref{cor2}}}

On suppose ici que les degrés dynamiques de $f$ vérifient $d_1 \leq
\dots \leq d_{s-1} < d_s > d_{s+1} \geq \dots \geq d_k$. Soit $\mu$ une mesure invariante, ergodique telle que $\log d(x
,\Acal) \in L^1(\mu)$ et $h(\mu) > \max ( \log d_{s-1}, \log
d_{s+1})$ (ou $h(\mu) > \log d_{k-1}$ si $s=k$). Montrons par l'absurde que $\chi_s >0$ et $\chi_{s+1} <0$.

Si $\chi_s \leq 0$ alors $\chi_{s-l}^{+}= \dots = \chi_s^{+}= \dots =
\chi_k^{+}=0$ et la première formule du théorème donnerait $\log
d_{s-1} < h(\mu) \leq \log d_{s-l-1}$ qui est absurde. De même si $\chi_{s+1} \geq 0$,
on applique la deuxième formule avec $s=s+1$ et on obtient $\log
d_{s+1} < h(\mu) \leq
\log d_{s+1+l'}$ qui est une contradiction (ici $\chi_1^{-}=
\dots = \chi_{s+1}^{-}= \dots = \chi_{s+1+l'}^{-}=0$). 

Passons à la minoration de $\chi_s$. Par la première formule du théorème et on a:
$$h(\mu) \leq \log d_{s-l-1} + 2(l+1) \chi_s$$
car par ce que l'on a fait précédemment on a $\chi_{s+1}^{+}= \dots = \chi_k^{+} =0$ et $\chi_s^{+}=
\chi_s$. On obtient $\chi_s \geq \frac{1}{2} \left(
\frac{1}{l+1} h(\mu)- \frac{1}{l+1} \log d_{s-l-1} \right)$. La concavité
de la fonction $q \rightarrow \log d_q$ (voir \cite{Gu} et \cite{Gr1}) implique $ \log d_{s-1} \geq \frac{1}{l+1} \log d_{s-l-1} + (1-
\frac{1}{l+1}) \log d_s$. On a donc
$$\chi_s \geq \frac{1}{2}
\left(\frac{1}{l+1} h(\mu) - \frac{1}{l+1} \log d_s - \log d_{s-1} +
\log d_s \right).$$
Mais comme $h(\mu) \leq \log d_{s}$ (voir \cite{DS1} et \cite{DS2}), cette dernière quantité
est supérieure à $\frac{1}{2} ( h(\mu) - \log d_s - \log d_{s-1} +
\log d_s)$, ce qui nous donne la minoration de $\chi_s$ que
l'on cherche.

Pour la majoration de $\chi_{s+1}$ la méthode est exactement la même à
condition d'utiliser la deuxième formule avec $s=s+1$.

\section{{\bf Théorie de Pesin et applications}}{\label{Pesin}}

Dans ce paragraphe, on considère une mesure de probabilité invariante,
ergodique $\mu$ telle que $\log d(x, \Acal) \in L^1(\mu)$ (avec $\Acal=C_f \cup I_f$). On va voir que cette
hypothèse permet de définir les exposants de Lyapounov pour $\mu$ et
de faire de la théorie de Pesin.

Tout d'abord, on définit l'extension naturelle $\widehat{X}$ de $X$
par:
$$ \widehat{X}:= \{ \widehat{x}=( \dots, x_{-n} , \dots, x_0, \dots ,
x_n , \dots) \in
X^{\Zz} \mbox{ , } f(x_{-n})=x_{-n+1} \}.$$
C'est l'ensemble des histoires des points de $X$. Dans cet espace
$f$ induit une application $\widehat{f}$ qui est le décalage à droite
et si on note $\pi$ la projection canonique $\pi(\widehat{x})=x_0$,
alors $\mu$ se relève en une unique probabilité $\widehat{\mu}$
invariante par $\widehat{f}$ qui vérifie $\pi_{*}
\widehat{\mu}= \mu$.

Dans l'espace $\widehat{X}$, on ne gardera que les orbites qui ne
visitent pas l'ensemble $\Acal$. On considère donc:

$$\widehat{X}^{*}= \{ \widehat{x} \in \widehat{X} \mbox{ , } x_n
\notin \Acal \mbox{ , } \forall n \in \Zz \}.$$

Cet ensemble est invariant par $\widehat{f}$ et $\widehat{\mu}(\widehat{X}^{*})=1$ car $\mu(\Acal)=0$.

Maintenant, on peut munir $X$ d'une famille de cartes $(\tau_x)_{x \in
  X}$ telles que $\tau_x(0)=x$, $\tau_x$ est définie sur une boule
  $B(0, \epsilon_0)$ avec $\epsilon_0$ indépendant de $x$ et la norme
  de la dérivée première et seconde de $\tau_x$ sur $B(0, \epsilon_0)$
  est majorée par une constante indépendante de $x$. Pour construire
  ces cartes il suffit de partir d'une famille finie $(U_i, \psi_i)$ de
  cartes de $X$ et de les composer par des translations.

Dans toute la suite, on notera $f_x= \tau_{f(x)}^{-1} \circ f \circ
\tau_x$ qui est définie au voisinage de $0$ quand $x$ n'est pas dans
$I_f$ et on posera aussi:
$$f_x^n= f_{f^{n-1}(x)} \circ \dots \circ f_x$$
$$f_{\widehat{x}}^{-n}= f_{x_{-n}}^{-1} \circ \dots \circ f_{x_{-1}}^{-1}$$
où dans $f_{x_{-i}}^{-1}$ on a pris la bonne branche inverse de $f$
par rapport à $\widehat{x}$:
celle qui envoie $x_{-i+1}$ sur $x_{-i}$.

Pour $\widehat{x} \in \widehat{X}^{*}$ on définit $D \widehat{f} (
\widehat{x})=Df_x(0)$ (où $\pi( \widehat{x})=x$). L'application $D
\widehat{f}$ va de $\widehat{X}^{*}$ dans $GL(k, \Cc)$ et c'est à ce
cocycle que nous allons appliquer la théorie de Pesin. Tout d'abord,
nous avons le:

\begin{lemme}

Les fonctions $\log^{+} \| D \widehat{f} ( \widehat{x}) \|$ et
$\log^{+} \|( D \widehat{f} ( \widehat{x}))^{-1} \|$ sont dans
$L^1(\widehat{\mu})$.

\end{lemme}

\begin{proof}

Commençons par montrer que $\log^{+} \| D \widehat{f} ( \widehat{x})
\|$ est dans $L^1(\widehat{\mu})$.

Si on applique le lemme 2.1 de l'article \cite{DiDu} de
T.-C. Dinh et C. Dupont à $f$, on obtient:

Il existe $\tau > 0$ et $p \in \Nn^{*}$ tels que pour tout $x$ hors de $I_f$:
$$ \| Df(x) \| + \| D^2 f (x) \| \leq \tau d(x, I_f)^{-p}.$$

On en déduit:
$$ \log \| Df_x(0) \| \leq \log \tau + \log d(x, \Acal)^{-p}$$
d'où
$$ \log^{+} \|Df_x(0) \| \leq \log \tau +\log d(x, \Acal)^{-p}$$
car on peut supposer $\tau > 1$ et que le diamètre de $X$ est
inférieur à $1$.
Mais
$$\int \log^{+} \| D \widehat{f} ( \widehat{x}) \| d
\widehat{\mu}(\widehat{x})= \int \log^{+} \| D f_{\pi(\widehat{x})}
(0) \| d \widehat{\mu}(\widehat{x})$$
qui est égal à 
$$\int \log^{+} \| D f_{x}(0) \| d (\pi_{*} \widehat{\mu})(\widehat{x})=\int \log^{+} \| D
f_{x}(0) \| d \mu(x).$$
Comme par hypothèse sur $\mu$, la fonction $\log d(x, \Acal)$ est dans
$L^1( \mu)$, on en déduit que la dernière intégrale ci-dessus est
finie. Autrement dit, la fonction $\log^{+} \| D \widehat{f} ( \widehat{x}) \|$ est bien dans $L^1( \widehat{\mu})$.

Passons maintenant à $\log^{+} \|( D \widehat{f} ( \widehat{x}))^{-1}
 \|$. Comme dans le lemme 2.1 de \cite{DiDu}, nous allons utiliser
 l'inégalité de Lojasiewicz (voir \cite{Lo}, IV.7.2 ). En effet cette inégalité nous
 donne:

Il existe $\tau' > 0$ et $p' \in \Nn^{*}$ tels que pour tout $x$ hors de $C_f$:
$$ | \mbox{Jac} f(x) |^2 \geq \tau' d(x, C_f)^{p'}.$$
Ensuite, grâce à l'estimée du lemme 2.1 de \cite{DiDu} donnée plus
haut, on sait que les modules des valeurs propres de $Df(x)^H Df(x)$ sont majorés
par $\tau^2 d(x,I_f)^{-2p}$ (ici la matrice $Df(x)^H$ est la
transposée-conjuguée de $Df(x)$). Le minimum des modules des valeurs propres de la
matrice $Df(x)^H Df(x)$ est donc minoré par
$\frac{\tau'}{\tau^{2(k-1)}} d(x,C_f)^{p'} d(x,I_f)^{2p(k-1)}$. Cela
implique que $\| (D f(x))^{-1} \|$ est majorée par $\tau'' d(x, C_f \cup
I_f)^{-p''}$ pour certains $\tau'' > 0$ et $p'' \in \Nn^{*}$.

On a donc:
$$\log \| (Df_x(0))^{-1} \| \leq \log \tau'' + \log d(x, \Acal)^{-p''}$$

qui est une fonction intégrable pour $\mu$.

Comme précédemment on en déduit que  $\log^{+} \|( D \widehat{f} (
\widehat{x}))^{-1} \| $ est dans $L^1( \widehat{\mu})$.

\end{proof}

Grâce à ce lemme on peut appliquer le théorème d'Oseledec et le
théorie de Pesin au cocycle $D \widehat{f}$. On a (voir \cite{KH}):

\begin{theoreme*} (Oseledec)

Il existe un ensemble invariant $\widehat{Y} \subset \widehat{X}^{*}$
avec $\widehat{\mu}(\widehat{Y})=1$ tel que pour tout $\widehat{x} \in
\widehat{Y}$ on ait:

1) L'existence d'une décomposition de $\Cc^k$:
$$\Cc^k= \oplus_{i=1}^{q} E_i(\widehat{x}),$$
où les $E_i(\widehat{x})$ vérifient $D \widehat{f}(\widehat{x})
E_i(\widehat{x})=E_i( \widehat{f}(\widehat{x}))$.

2) L'existence de fonctions (exposants de Lyapounov de $f$):
$$\lambda_1 > \dots > \lambda_q,$$
avec
$$\lim_{m \rightarrow \pm \infty} \frac{1}{|m|} \log \|
  A(\widehat{x},m)v \|= \pm \lambda_i$$
pour tout $v \in E_i(\widehat{x}) \setminus \{ 0 \}$.
Ici les $A(\widehat{x},m)$ sont definis par:
$$A(\widehat{x},m)= D \widehat{f}(\widehat{f}^{m-1}(\widehat{x})) \circ
\dots \circ D \widehat{f}(\widehat{x}) \mbox{ pour } m>0$$
$$A(\widehat{x},0)=Id$$
$$A(\widehat{x},m)=D \widehat{f}(\widehat{f}^{m}(\widehat{x}))^{-1}
\circ \dots \circ D \widehat{f}(\widehat{f}^{-1}(\widehat{x}))^{-1} \mbox{ pour } m<0.$$

\end{theoreme*}

De plus, nous avons le:

\begin{theoreme*} ($\gamma$-réduction de Pesin)

Pour tout $\gamma >0$ il existe une fonction $C_{\gamma}: \widehat{X}^{*}
\rightarrow GL(k, \Cc)$ telle que:

1) $\lim_{m \rightarrow \infty} \frac{1}{m} \log \| C_{\gamma}^{\pm 1
   }(\widehat{f}^m(\widehat{x})) \|=0$ (on parle de fonction
   tempérée).

2) Pour presque tout $\widehat{x}$, la matrice
   $A_{\gamma}(\widehat{x})=C_{\gamma}^{-1}(\widehat{f}(\widehat{x})) D
   \widehat{f}( \widehat{x}) C_{\gamma}(\widehat{x})$ a la forme
   suivante:

$$A_{\gamma}(\widehat{x})= \left(
\begin{array}{cccc}
A_{\gamma}^1(\widehat{x}) \\
 & \ddots & \\
 &  & A_{\gamma}^{q}(\widehat{x})\\
\end{array}
\right)$$
 où chaque $A_{\gamma}^i(\widehat{x})$ est une matrice carrée de
 taille $dim E_i(\widehat{x}) \times dim E_i(\widehat{x})$ et on a:
$$e^{\lambda_i- \gamma} \leq \| A_{\gamma}^i(
 \widehat{x})^{-1} \| ^{-1} \mbox{  et  } \| A_{\gamma}^i(
 \widehat{x}) \| \leq e^{\lambda_i+ \gamma}.$$

3) Enfin, pour presque tout $\widehat{x}$ l'application $C_{\gamma}(\widehat{x})$ envoie la
   décomposition standard $\oplus_{i=1}^{q}
   \Cc^{dim E_i(\widehat{x})}$ sur $\oplus_{i=1}^{q} E_i(\widehat{x})$.

\end{theoreme*} 

Dans toute la suite nous noterons $\widehat{Y}$ l'ensemble des points
de $\widehat{X}^{*}$ qui vérifient les conclusions de ces deux
théorèmes. De plus, pour $\widehat{x} \in \widehat{Y}$, nous
appellerons $\chi_1 \geq \dots \geq \chi_k$
les exposants de Lyapounov de $\widehat{\mu}$ notés avec répétition
(contrairement aux $\lambda_i$ du théorème précédent).

Notons maintenant $g_{\widehat{x}}$ la lecture de $f_x$ dans ces
cartes (i.e. $g_{\widehat{x}}=C_{\gamma}^{-1}(
\widehat{f}(\widehat{x})) \circ f_x \circ C_{\gamma}(\widehat{x})$
avec $\pi(\widehat{x})=x$).

Nous allons donner quelques propriétés de $g_{\widehat{x}}$
(avec $\widehat{x} \in \widehat{Y}$) qui nous seront utiles pour la
suite.

\begin{proposition}{\label{direct}}

Il existe des constantes $\tau > 0$, $\epsilon_0 >0$ et $p \in
\Nn^{*}$ qui ne dépendent que de $f$ et $X$ telles que:

1) $g_{\widehat{x}}(0)=0$

2) $Dg_{\widehat{x}}(0)= \left(
\begin{array}{cccc}
A_{\gamma}^1(\widehat{x}) \\
 & \ddots & \\
 &  & A_{\gamma}^{q}(\widehat{x})\\
\end{array}
\right)$
 
3) Si on note $g_{\widehat{x}}(w)=Dg_{\widehat{x}}(0)w + h(w)$, on a:
$$\| Dh(w) \| \leq \tau \| C_{\gamma}( \widehat{f}(\widehat{x}))^{-1} \| \|
C_{\gamma}(\widehat{x}) \|^{2} d(x, \Acal)^{-p} \| w \|$$
pour $\|w\| \leq \epsilon_0 d(x, \Acal)/ \|C_{\gamma}(\widehat{x})\|$.

\end{proposition}

\begin{proof}

Le premier point est évident. Le second provient immédiatement du
théorème précédent. Il reste à prouver la troisième propriété.

On a $Dg_{\widehat{x}}(w)=Dg_{\widehat{x}}(0) + Dh(w)$, d'où:
$$\| Dh(w) \| = \| Dg_{\widehat{x}}(w) - Dg_{\widehat{x}}(0) \| \leq
\| C_{\gamma}( \widehat{f}(\widehat{x}))^{-1} \| \| Df_x
(C_{\gamma}(\widehat{x})w) -Df_x(0) \|  \|C_{\gamma}(\widehat{x})
\|.$$
Par l'estimée sur $\| D^2 f \|$ de T.-C. Dinh et C. Dupont qui se trouve au début de la
démonstration du lemme précédent, on a:

$$ \| Df_x (C_{\gamma}(\widehat{x})w) -Df_x(0) \| \leq \tau d( \tau_x
([0, C_{\gamma}(\widehat{x})w]), \Acal)^{-p} \|
C_{\gamma}(\widehat{x})w \|$$
pour $\|C_{\gamma}(\widehat{x})w \| \leq \epsilon_0$.
L'image par $\tau_x$ du segment $[0, C_{\gamma}(\widehat{x})w]$ vit
dans la boule $B(x, K \|C_{\gamma}(\widehat{x})w \|)$ où $K$ est une
constante qui ne dépend que de $X$. On en déduit que pour
$\|C_{\gamma}(\widehat{x})w \| \leq \epsilon_0 d(x, \Acal)$,
on a:
$$\| Dh(w) \|  \leq \| C_{\gamma}( \widehat{f}(\widehat{x}))^{-1} \|
\|C_{\gamma}(\widehat{x}) \|^2 \tau (1-K \epsilon_0)^{-p} d(x, \Acal)^{-p}
\|w \|.$$

Quitte à prendre $\epsilon_0$ petit, on peut supposer que $K
\epsilon_0$ est inférieur à $\frac{1}{2}$ et donc modulo un
changement de la constante $\tau$ on a la majoration voulue de $\| Dh(w) \|$.

\end{proof}

Dans la démonstration des formules du théorème nous utiliserons aussi
$g_{\widehat{x}}^{-1}=C_{\gamma}^{-1}(\widehat{f}^{-1}(\widehat{x}))
  \circ f_{\widehat{x}}^{-1} \circ C_{\gamma}(\widehat{x})=C_{\gamma}^{-1}(\widehat{f}^{-1}(\widehat{x}))
  \circ f_{x_{-1}}^{-1} \circ C_{\gamma}(\widehat{x})$ et des
  estimées sur cette application qui sont données par la:

\begin{proposition}{\label{inverse}}

Il existe des constantes $\tau > 0$, $\epsilon_0 >0$ et $p \in
\Nn^{*}$ qui ne dépendent que de $f$ et $X$ telles que:

1) $g_{\widehat{x}}^{-1}(w)$ est
bien définie pour $\|w\| \leq \epsilon_0 d(x_{-1}, \Acal)^p / \|C_{\gamma}(\widehat{x})\|$.

2) $g_{\widehat{x}}^{-1}(0)=0$

3) $Dg_{\widehat{x}}^{-1}(0)= \left(
\begin{array}{cccc}
(A_{\gamma}^1(\widehat{f}^{-1}(\widehat{x})))^{-1} \\
 & \ddots & \\
 &  & (A_{\gamma}^{q}(\widehat{f}^{-1}(\widehat{x})))^{-1}\\
\end{array}
\right)$
 
4) Si on note $g_{\widehat{x}}^{-1}(w)=Dg_{\widehat{x}}^{-1}(0)w + h(w)$, on a:
$$\| Dh(w) \| \leq \tau \| C_{\gamma}( \widehat{f}^{-1}(\widehat{x}))^{-1} \| \|
C_{\gamma}(\widehat{x}) \|^{2} d(x_{-1}, \Acal)^{-p} \| w \|$$
pour $\|w\| \leq \epsilon_0 d(x_{-1}, \Acal)^p / \|C_{\gamma}(\widehat{x})\|$.

\end{proposition}

\begin{proof}

Commençons par démontrer que $g_{\widehat{x}}^{-1}(w)$ est bien défini
pour 
$$\|w\| \leq \epsilon_0 d(x_{-1}, \Acal)^p / \|C_{\gamma}(\widehat{x})\|.$$

Cela va reposer sur le lemme 2 de \cite{BD1} (construction de branches inverses).

Rappelons que grâce à l'estimée de T.-C. Dinh et C. Dupont qui se trouve au début de la
démonstration du lemme précédent on a l'existence de $\tau' > 0$ et $p' \in \Nn^{*}$ tels que pour tout $x$ hors de $I_f$:
$$ \| Df(x) \| + \| D^2 f (x) \|  \leq \tau' d(x, I_f)^{-p'}$$

D'autre part, dans la preuve de ce lemme, on a obtenu aussi l'existence de $\tau'' >0$ et $p'' \in \Nn^{*}$ avec $\| (D f(x))^{-1} \|
\leq \tau'' d(x, C_f \cup I_f)^{-p''}$ pour $x$ hors de $\Acal=C_f \cup I_f$.

Quitte à remplacer $p'$ et $p''$ par le maximum des deux, on pourra
supposer dans la suite que $p'=p''$ et de même $\tau'=\tau''$. Ce sont
des constantes qui ne dépendent que de $f$ et $X$.

Soit $w$ tel que $\| w \| \leq \epsilon_0' d(x_{-1}, \Acal)$. On a
$\tau_{x_{-1}}(w) \in B(x_{-1}, K \epsilon_0' d(x_{-1}, \Acal))$ (où $K$ ne dépend que de $X$) et alors:
$$ \| D f_{x_{-1}}(w) \| + \| D^2 f_{x_{-1}}(w) \| + \| (D
f_{x_{-1}}(w))^{-1} \| \leq  2 \tau' (1-K \epsilon_0')^{-p'} d(x_{-1}, \Acal)^{-p'}\leq \tau' d(x_{-1}, \Acal)^{-p'},$$
si on prend $\epsilon_0'$ petit et quitte à renommer $\tau'$.

L'inégalité ci-dessus combinée avec le lemme 2 de \cite{BD1} implique que
$f_{\widehat{x}}^{-1}$ est définie sur une boule $B(0, \epsilon_0'
d(x_{-1},\Acal)^{3p'}/ {\tau'}^3)=B(0,2 \epsilon_0 d(x_{-1}, \Acal)^p)$. En particulier, $g_{\widehat{x}}^{-1}(w)=C_{\gamma}^{-1}(\widehat{f}^{-1}(\widehat{x}))
  \circ f_{\widehat{x}}^{-1} \circ C_{\gamma}(\widehat{x})(w)$ est
bien défini pour $\|w\| \leq 2 \epsilon_0 d(x_{-1},
\Acal)^{p} / \|C_{\gamma}(\widehat{x})\|$.

Passons maintenant à la majoration de $\| D h(w) \|$. Pour cela, comme dans
la proposition précédente, il faut contrôler $\| D^2
f_{\widehat{x}}^{-1} \|$.

L'image de $B(0,  2 \epsilon_0 d(x_{-1},\Acal)^{p})$ par
$f_{\widehat{x}}^{-1}$ est incluse dans $B(0,1)$ (toujours par le lemme 2 de \cite{BD1}). La
formule de Cauchy nous donne donc:

$$ \| D f_{\widehat{x}}^{-1}(w) \| + \| D^2 f_{\widehat{x}}^{-1}(w) \|
\leq \tau'' d(x_{-1}, \Acal)^{-p''},$$
pour $\|w\| \leq \epsilon_0 d(x_{-1}, \Acal)^{p}$.

Grâce à cette majoration, on peut maintenant contrôler $\| D h(w)
\|$. On a:
$$ \| Dh(w) \| = \| Dg_{\widehat{x}}^{-1}(w)- Dg_{\widehat{x}}^{-1}(0)
\| \leq \| C_{\gamma}( \widehat{f}^{-1}(\widehat{x}))^{-1} \| \|
Df_{\widehat{x}}^{-1}(C_{\gamma}(\widehat{x})w)-
Df_{\widehat{x}}^{-1}(0)\| \|C_{\gamma}(\widehat{x})\|, $$

d'où:

$$\| Dh(w) \| \leq \tau'' \| C_{\gamma}( \widehat{f}^{-1}(\widehat{x}))^{-1} \| \|
C_{\gamma}(\widehat{x}) \|^{2} d(x_{-1}, \Acal)^{-p''} \| w \|,$$
pour $\|w\| \leq \epsilon_0 d(x_{-1}, \Acal)^p / \|C_{\gamma}(\widehat{x})\|.$

Cela démontre bien la proposition quitte à prendre le maximum entre
$p$ et $p''$.

\end{proof}

Dans les deux propositions, on voit que la distance de $x$ à
l'ensemble $\Acal$ joue un rôle crucial. Comme on les utilisera le
long d'orbites de point, on aura besoin de savoir la distance entre
$f^i(x)$ et $\Acal$. Celle-ci est donnée par le:

\begin{lemme}{\label{distance}}

Il existe un ensemble $\widehat{Y}$ dans $\widehat{X}$ de mesure
pleine pour $\widehat{\mu}$ tel que pour tout $\widehat{x} \in
\widehat{Y}$ on ait:
$$d(x_n, \Acal) \geq V(\widehat{x}) e^{- \gamma |n|}$$
pour tout $n \in \Zz$. Ici $V$ est une fonction mesurable à valeur dans $\Rr^{*}_{+}$.

\end{lemme}

\begin{proof}

Il suffit d'appliquer le théorème de Birkhoff à la fonction
$u(\widehat{x}) =\log d(\pi(\widehat{x}), \Acal)$ qui est dans
$L^1(\widehat{\mu})$ (voir par exemple le Lemme 2.3 de \cite{DiDu}).

\end{proof}

Dans toute la suite, $\widehat{Y}$ désignera le sous-ensemble de
points de $\widehat{X}^{*}$ qui vérifient les théorèmes d'Oseledec, de
$\gamma$-réduction de Pesin et les conclusions du lemme précédent.

Nous allons maintenant faire des rappels sur la transformée de graphe.

\section{{\bf Transformée de graphe}}{\label{graphe}}

La cadre de ce paragraphe est $\Cc^k$. Dans toute la suite $\|.\|$
désignera la norme Euclidienne.

On considère 
$$g(X,Y)=(AX + R(X,Y), BY + U(X,Y))$$
avec $(X,0) \in E_1$ (abscisses), $(0,Y) \in E_2$ (ordonnées) et $A$, $B$ des matrices. On suppose aussi que $g(0,0)=(0,0)$ et $\max (\| D
R(Z) \|, \| D U (Z) \|) \leq \delta$ dans la boule
$B(0,r)$. Enfin par hypothèse, on aura $\gamma \leq \| A \| \leq \| B^{-1}
\|^{-1}(1 - \gamma)$. Soit maintenant $\{ ( \Phi(Y),Y), Y \in D \}$ un
graphe dans $B(0,r)$ au-dessus d'une partie $D$ de $E_2$ qui vérifie $\| \Phi(Y_1) - \Phi(Y_2) \| \leq \gamma_0 \| Y_1 - Y_2
\|$. Dans le théorème qui suit, on donne des conditions sur $\delta$,
$\gamma$ et $\gamma_0$ pour que l'image de ce graphe par $g$ soit un
graphe qui vérifie le même contrôle.

\begin{theoreme*}

Si $\delta \| B^{-1} \| (1 + \gamma_0) <1$ alors l'image par $g$ du
graphe précédent est un graphe au-dessus de $\pi_0(g(\mbox{graphe de }
\Phi))$ où $\pi_0$ est la projection sur les ordonnées. Par
ailleurs, si $(\Psi(Y),Y)$ désigne ce nouveau graphe, on a:
$$\| \Psi(Y_1) - \Psi(Y_2) \| \leq \frac{\| A \| \gamma_0 + \delta(1+
  \gamma_0)}{\| B^{-1} \|^{-1} - \delta(1+ \gamma_0)}  \| Y_1 - Y_2
  \|$$
qui est inférieur à $\gamma_0 \| Y_1 - Y_2 \|$ si $\delta \leq
  \epsilon(\gamma_0, \gamma)$.

Enfin, si de plus $B(0 , \alpha) \subset D$ et $\| \Phi(0) \|
\leq \beta$, alors $\pi_0(g(\mbox{graphe de } \Phi))$ contient $B(0,
(\| B^{-1} \|^{-1} - \delta(1 + \gamma_0)) \alpha - \| A \| \beta -
\delta \beta - \|
D^2 g \|_{B(0,r)} \beta^2)$ et $\| \Psi(0) \| \leq   (1+ \gamma_0)( \| A \| \beta + \delta \beta +
\|D^2g \|_{B(0,r)} \beta^2) $ (si $\delta \leq \epsilon(\gamma_0, \gamma)$).

\end{theoreme*}

\begin{proof}

La démonstration est tirée essentiellement de \cite{KH} mais nous
préférons la donner par confort pour le lecteur.
 
Soit $\lambda(Y)= BY + U(\Phi(Y),Y)$. C'est l'ordonnée de
$g(\Phi(Y),Y)$. Pour démontrer que $g(\mbox{graphe de } \Phi)$ est un
graphe au-dessus de $\pi_0(g(\mbox{graphe de } \Phi))$, il suffit de
voir que $\lambda$ est une bijection de $D$ sur $\lambda(D)$ c'est-à-dire que $\lambda(Y)=Y_0$ a une unique solution pour $Y_0 \in
\lambda(D)$.

Posons $\gamma(Y)=B^{-1} Y_0 - B^{-1} U(\Phi(Y),Y)$. On a que
$\lambda(Y)=Y_0$ est équivalent à $\gamma(Y)=Y$. Mais, si $Y_1 \mbox{
  , } Y_2$ sont dans $D$, on a:
$$\| \gamma(Y_1) - \gamma(Y_2) \| \leq  \| B^{-1} \| \|
U(\Phi(Y_1),Y_1)- U(\Phi(Y_2),Y_2) \|$$
qui est inférieur à $ \| B^{-1}  \| \delta (1 + \gamma_0) \|Y_1 -
Y_2 \|$. Autrement dit, quand $\delta \| B^{-1} \| (1 + \gamma_0) <1$,
alors $\gamma(Y)=Y$ a bien une unique solution dans $D$.

Passons maintenant au contrôle de la pente du graphe $(\Psi(Y),Y)$ que
l'on a obtenu.

On considère $(\Psi(Y_1'),Y_1')$ et $(\Psi(Y_2'),Y_2')$ deux points du
graphe. Ils sont l'image de $(\Phi(Y_1),Y_1)$ et $(\Phi(Y_2),Y_2)$ par
$g$. On notera $X_i'= \Psi(Y_i')$ ($i=1,2$). Alors, d'une part:
$$\| X_1' - X_2' \| = \| A \Phi(Y_1) + R(\Phi(Y_1),Y_1) - A
\Phi(Y_2) - R(\Phi(Y_2),Y_2)\|$$
qui est inférieur à
$$( \| A \| \gamma_0 + \delta(1 + \gamma_0)) \|Y_1 - Y_2 \|.$$
D'autre part:

$$ \| Y_1' - Y_2' \| = \| B Y_1 + U(\Phi(Y_1),Y_1) - B Y_2 -
U(\Phi(Y_2),Y_2) \|$$
qui est supérieur à
$$( \| B^{-1} \|^{-1} - \delta (1 + \gamma_0)) \| Y_1 - Y_2 \|$$
car $\| B(Y_1-Y_2) \| \geq  \| B^{-1} \|^{-1} \| Y_1 - Y_2 \|$.

En combinant ces deux inégalités, on obtient:

$$\| \Psi(Y_1') - \Psi(Y_2') \| \leq \frac{\| A \| \gamma_0 + \delta(1
  + \gamma_0)}{ \| B^{-1} \|^{-1} - \delta (1 + \gamma_0)}\| Y_1' -
Y_2' \|$$
qui est l'inégalité cherchée.

Passons maintenant aux dernières estimations.

Tout d'abord on veut majorer la distance entre $(0,0)$ et l'image de
$(\Phi(0),0)$ par $g$. Soit $v$ le vecteur $(\Phi(0),0)$ normalisé de
sorte que $\| v \| = 1$. Si $P$ est un point du segment $[(0,0),
  (\Phi(0),0)]$, on a:
$$\| Dg(P)v - Dg(0)v \| \leq \|D^2g \|_{B(0,r)} \|P\|$$
Mais $\| Dg(0)v \| \leq \| A \| + \delta$ donc:
$$\| Dg(P)v \| \leq  \| A \| + \delta +  \|D^2g \|_{B(0,r)} \beta.$$
On déduit de cette inégalité que la distance entre $(0,0)$ et l'image de
$(\Phi(0),0)$ par $g$ est majorée par $ \| A \| \beta + \delta \beta +
\|D^2g \|_{B(0,r)} \beta^2$.

Maintenant, on a:
$$ \| \lambda(Y_1) - \lambda(Y_2) \| \geq \| B(Y_1 - Y_2) \| - \delta
(1+\gamma_0) \|Y_1 -Y_2 \|$$
qui est plus grand que $(\| B ^{-1} \|^{-1} - \delta(1 + \gamma_0)) \|
Y_1 - Y_2 \|$. En particulier, $\| \lambda(Y) - \lambda(0) \| \geq (\|
B ^{-1} \|^{-1} - \delta(1 + \gamma_0))  \alpha$ pour $Y \in
\partial B(0, \alpha)$. Comme $\| \lambda(0) \|$ est inférieur à $ \| A \| \beta + \delta \beta +
\|D^2g \|_{B(0,r)} \beta^2$, on en déduit que $\pi_0(g(\mbox{graphe de } \Phi))$ contient $B(0,
(\| B^{-1} \|^{-1} - \delta(1 + \gamma_0)) \alpha - \| A \| \beta -
\delta \beta - \| D^2 g \|_{B(0,r)} \beta^2)=B(0,r')$.

Il reste à majorer $\| \Psi(0) \|$. On a:
$$ \| \Psi(0) - \Psi(\lambda(0)) \| \leq  \gamma_0 \| \lambda(0) \|$$
(si on suppose que $\delta \leq \epsilon(\gamma_0, \gamma)$).
On obtient alors (toujours avec la majoration de la distance entre
$(0,0)$ et $(\Psi(\lambda(0)), \lambda(0))=g(\Phi(0),0)$)
$$ \| \Psi(0) \| \leq  (1 +\gamma_0)( \| A \| \beta + \delta \beta +
\|D^2g \|_{B(0,r)} \beta^2).$$
C'est l'estimée que l'on cherchait.

\end{proof}

On va passer maintenant à la démonstration des deux formules.

\section{{\bf Démonstration de la première inégalité du théorème}}

Commençons par rappeler la définition de l'entropie métrique. 

Notons $d_n(x,y)= \max_{0  \leq i \leq n-1} \{ d(f^i(x), f^i(y)) \}$ et $B_n(x, \delta)$ la
boule de centre $x$ et de rayon $\delta$ pour cette métrique. Par le
théorème de Brin et Katok (voir \cite{BK}), l'entropie
métrique de $\mu$ est donnée par la formule:
$$h(\mu)= \lim_{\delta \rightarrow 0} \liminf_{n} - \frac{1}{n} \log \mu(B_n(x,
\delta))$$
pour $\mu$-presque tout $x$.

On va maintenant faire quelques uniformisations. 

Soit $\Lambda_{\delta, n}= \{ x \mbox{ , } \mu (B_n(x,\delta)) \leq
e^{-    h(\mu)n + \gamma n} \}$.

Si $\delta$ est pris petit on a 
$$\frac{4}{5} \leq \mu( \{ x \mbox{ , }  \liminf_{n} - \frac{1}{n} \log \mu(B_n(x,\delta)) \geq
  h(\mu) - \frac{\gamma}{2} \}) \leq \mu( \cup_{n_0} \cap_{n \geq n_0}
  \Lambda_{\delta, n }).$$
En particulier, si on prend $n_0$ grand, on a: $\mu(\cap_{n \geq n_0} \Lambda_{\delta, n}) \geq 3/4$.

Rappelons que l'on note $\widehat{Y}$ l'ensemble des bons points de
$\widehat{X}^{*}$ pour la théorie de Pesin.

Posons 
$$\widehat{Y}_{\alpha_0}=\{  \widehat{x} \in \widehat{Y} \mbox{
  , } \alpha_0 \leq \| C_{\gamma} ( \widehat{x})^{\pm 1} \| \leq
  \frac{1}{\alpha_0} \mbox{ , } V(\widehat{x}) \geq \alpha_0 \}$$
(voir le paragraphe \ref{Pesin} pour les notations).
Si $\alpha_0$ est suffisamment petit, on a
  $\widehat{\mu}(\widehat{Y}_{\alpha_0}) \geq 3/4$ d'où $\mu(A_{n_0}) \geq
  1/2$ avec $A_{n_0}=  \pi(\widehat{Y}_{\alpha_0}) \cap (
  \cap_{n \geq n_0} \Lambda_{\delta,n})$.

Maintenant, comme pour les points $x$ de $A_{n_0}$ on a $\mu( B_n(x,
\delta)) \leq e^{-h(\mu)n + \gamma n}$, on peut trouver un ensemble
$\{x_i \}_{1 \leq i \leq N}$ d'éléments de $A_{n_0}$ avec $N \geq
\frac{1}{2}e^{h(\mu)n - \gamma n }$ tels que $x_i =
\pi(\widehat{x_i})$ où $\widehat{x_i} \in \widehat{Y}_{\alpha_0}$ et
avec les $B_n(x_i, \delta/2)$ disjointes (i.e. les points $x_i$ sont $(n,\delta)$-séparés).

$ $

Voici le plan de la démonstration de la formule. Dans celle-ci, on
adapte des idées de J. Buzzi (voir \cite{Bu}) et S. E. Newhouse (voir
\cite{Ne}) à notre contexte: celui des applications méromorphes. En chaque point $x_i$ nous allons
construire une variété stable approchée $W_i$ de dimension
$k-(s-l)+1$ (dans tout ce texte les dimensions seront des dimensions complexes). Cela signifiera en particulier que le diamètre de
$f^{q}(W_i)$ restera inférieur à $\delta/4$ ($q=0, \dots , n-1$) et que les $W_i$ seront
assez plates. Ensuite, dans un deuxième paragraphe nous minorerons le
volume $k-s+l+1$-dimensionnel de ces variétés par à peu près $e^{-2
  \chi_{s-l}^{+} n - \dots - 2 \chi_k^{+} n }$. Le volume total de
toutes ces variétés est donc essentiellement supérieur à
$e^{h(\mu)n-2  \chi_{s-l}^{+} n - \dots - 2 \chi_k^{+} n }$. Les $W_i$
étant assez plates, on pourra trouver un plan $P$ de dimension
$k-s+l+1$ tel que d'une part la projection de tous les $W_i$ sur $P$ soit de
volume minoré par la même quantité et d'autre part les $W_i$ seront des
graphes au-dessus de $P$. Maintenant, si $\pi_1$ désigne la
projection orthogonale sur $P$, la minoration de volume implique que les
fibres de $\pi_1$ (qui sont des plans de dimension $s-l-1$) coupent
$\cup W_i$ en moyenne en au moins $e^{h(\mu)n-2  \chi_{s-l}^{+} n -
  \dots - 2 \chi_k^{+} n }$ points. Mais vu que les $x_i$ sont $(n,
\delta)$-séparés et que les diamètres des poussés en avant des $W_i$
restent petits, l'intersection d'une fibre de $\pi_1$ avec $\cup W_i$
donne des points $(n, \delta/2)$-séparés. Cela signifie qu'en moyenne le
nombre de points $(n , \delta/2)$-séparés dans une fibre de $\pi_1$ est
minoré par $e^{h(\mu)n-2  \chi_{s-l}^{+} n -  \dots - 2 \chi_k^{+} n
}$. Enfin, dans le troisième paragraphe
nous donnerons une majoration de cette moyenne par essentiellement
$(\max_{0 \leq q \leq s-l-1} d_q)^n$ et cela prouvera l'inégalité.

\subsection{{\bf Construction des variétés stables approchées}}

Rappelons que l'on a $\chi_1 \geq \dots > \chi_{s-l}= \dots = \chi_s
\geq \dots \geq \chi_k$.
On notera $E_1(\widehat{x}) , \dots , E_{m}(\widehat{x})$ les
$E_i(\widehat{x})$ du théorème d'Oseledec correspondant aux exposants $\chi_1 , \dots , \chi_{s-l-1}$ et
$E_{m+1}(\widehat{x}) , \dots , E_{q}(\widehat{x})$ les
$E_i(\widehat{x})$ de $\chi_{s-l} , \dots , \chi_k$ (voir le
paragraphe \ref{Pesin} pour les notations). Soit:
$$E_u (\widehat{x}) = \oplus_{i=1}^{m} E_i(\widehat{x}) \mbox{ et }
E_s (\widehat{x}) = \oplus_{i=m+1}^{q} E_i(\widehat{x}).$$
Par ailleurs, $E_s (\widehat{x})$ sera dans la suite coupé en deux parties. Soit $n_1$ le nombre d'exposants parmi
$\chi_{s-l}, \dots , \chi_k$ qui sont strictement négatifs (bien sûr
$n_1$ peut être égal à $0$). Alors, nous noterons $E_s^1(
\widehat{x})$ la somme directe des $E_i(\widehat{x})$ ($i=m+1 , \dots ,
q$) correspondant aux $\chi_i$ strictement négatifs et
$E_s^2(\widehat{x})$ la somme directe des autres $E_i(\widehat{x})$ de
$E_s (\widehat{x})$. La dimension de $E_s^1(\widehat{x})$ est donc
$n_1$ et celle de $E_s^2(\widehat{x})$ est $k-s+l+1-n_1$.

Soit $x$ un des $N$ points $x_i$ ($x=\pi(\widehat{x})$ avec
$\widehat{x} \in \widehat{Y}_{\alpha_0}$). On va construire une
variété stable approchée qui passe par $x$ en utilisant la transformée
de graphe. Fixons $\gamma_0 >
0$ très petit devant $\alpha_0$. Dans toute la suite $n$ sera pris grand par rapport à des
constantes qui dépendent de $\gamma_0$ et $\gamma$.

On se place maintenant dans $C_{\gamma}^{-1}(\widehat{f}^n(\widehat{x}))  E_u (
\widehat{f}^n(\widehat{x})) \oplus
C_{\gamma}^{-1}(\widehat{f}^n(\widehat{x})) E_s ( \widehat{f}^n(\widehat{x}))$
et on part de 
$$ \{ 0 \}^{s-l-1}  \times B_2(0, e^{-4p \gamma n}) \times B_1(0, e^{-3p \gamma n}),$$
où $B_2(0, e^{-4p \gamma n})$ est la boule de $\Cc^{k-s+l+1-n_1}$ de centre $0$ et de
rayon $e^{-4p \gamma n}$ et $B_1(0, e^{-3p \gamma n})$ est celle de
$\Cc^{n_1}$ de centre $0$ et de rayon $e^{-3p \gamma n}$ . Cet ensemble est un graphe $(\Phi_n(Y),Y)$
au-dessus d'une partie de $C_{\gamma}^{-1}(\widehat{f}^n(\widehat{x}))
E_s ( \widehat{f}^n(\widehat{x}))$ (avec $\Phi_n(Y)=0$). 

\begin{lemme}

L'image du graphe $(\Phi_n(Y),Y)$ par $g_{\widehat{f}^n(
  \widehat{x})}^{-1}$ est un graphe $(\Phi_{n-1}(Y),Y)$ au-dessus
  d'une partie de $C_{\gamma}^{-1}(\widehat{f}^{n-1}(\widehat{x}))
E_s ( \widehat{f}^{n-1}(\widehat{x}))$. Il vérifie de plus $\|
  \Phi_{n-1}(Y_1) - \Phi_{n-1}(Y_2) \| \leq \gamma_0 \| Y_1 - Y_2 \|$.

\end{lemme}

\begin{proof}

Il s'agit d'utiliser le théorème du paragraphe \ref{graphe}.

Tout d'abord, si on prend pour abscisse $C_{\gamma}^{-1}(\widehat{f}^n(\widehat{x}))  E_u (
\widehat{f}^n(\widehat{x}))$ et pour ordonnée $C_{\gamma}^{-1}(\widehat{f}^n(\widehat{x}))  E_s (
\widehat{f}^n(\widehat{x}))$, on a 
$$g_{\widehat{f}^n( \widehat{x})}^{-1}(X,Y)=(AX + R(X,Y), BY +
U(X,Y))$$
avec
$$\gamma \leq \|A\| \leq \| (A_{\gamma}^{m} ( \widehat{f}^{n-1}(\widehat{x})))^{-1}
\| \leq (1 - \gamma) \| (A_{\gamma}^{m+1} (
\widehat{f}^{n-1}(\widehat{x}))) \|^{-1}=(1 - \gamma) \| B^{-1} \|^{-1}$$
par la proposition \ref{inverse}, le théorème de $\gamma$-réduction de
Pesin et le fait que $\gamma$ peut être supposé petit par rapport à
des constantes qui ne dépendent que des exposants de Lyapounov de $\mu$. 

De plus, toujours par cette proposition,
$$ \max(\| D R(X,Y) \|, \| D U(X,Y) \|) \leq \tau \|
C_{\gamma}(\widehat{f}^{n-1}(\widehat{x}))^{-1} \| \| C_{\gamma}(
\widehat{f}^{n}(\widehat{x})) \|^{2} d(f^{n-1}(x), \Acal)^{-p} \|(X,Y)\|$$
avec $\|(X,Y) \| \leq \frac{\epsilon_0 d(f^{n-1}(x), \Acal)^p}{ \|
  C_{\gamma}(\widehat{f}^n(\widehat{x}))\|}$.
Mais comme les fonctions $\|C_{\gamma}^{\pm 1} \|$ sont tempérées et que
$\widehat{x}$ est dans $\widehat{Y}_{\alpha_0}$, on peut supposer que $ \|
C_{\gamma}(\widehat{f}^{n-1}(\widehat{x}))^{-1} \| \| C_{\gamma}(
\widehat{f}^{n}(\widehat{x})) \|^{2} \leq \frac{1}{\alpha_0^3}
e^{3 \gamma n }$ (voir \cite{KH} p. 668) et on a $d(f^{n-1}(x) , \Acal)
\geq \alpha_0 e^{-n\gamma}$ (voir le lemme \ref{distance}).

Pour $\|(X,Y) \| \leq e^{-2p \gamma n}$, on obtient:
 
$$ \max(\| D R(X,Y) \|, \| D U(X,Y) \|) \leq e^{4 \gamma n}
e^{\gamma n p } e^{- 2 \gamma n p}$$
qui est très petit car $p$ peut être supposé supérieur à $5$. Cette quantité
joue le rôle de $\delta$ dans le théorème du paragraphe
\ref{graphe}. Comme  il est aussi petit que l'on veut pourvu que $n$ soit
grand, on a bien démontré le lemme.

\end{proof}

Maintenant, de ce graphe $(\Phi_{n-1}(Y),Y)$, on ne garde que la
partie qui se trouve au-dessus de $\{ 0 \}^{s-l-1} \times B_2(0, e^{- 3
  \gamma p n}) \times B_1(0, e^{- 3  \gamma p n}) $ (on fait un cut-off). Puis on prend son
image par $g_{\widehat{f}^{n-1}(\widehat{x})}^{-1}$ qui de nouveau est
un graphe $(\Phi_{n-2}(Y),Y)$ au-dessus d'une partie de $C_{\gamma}^{-1}(\widehat{f}^{n-2}(\widehat{x}))
E_s ( \widehat{f}^{n-2}(\widehat{x}))$ avec $\|
\Phi_{n-2}(Y_1) - \Phi_{n-2}(Y_2) \| \leq \gamma_0 \| Y_1 - Y_2 \|$
(la démonstration est exactement la même que dans le lemme précédent).

De ce graphe, on ne garde que la partie au-dessus $\{ 0 \}^{s-l-1} \times  B_2(0, e^{- 3
  \gamma p n}) \times B_1(0, e^{- 3  \gamma p n})$ et on continue ainsi le procédé jusqu'à
  obtenir un graphe $(\Phi_0(Y),Y)$ au-dessus d'une partie de $C_{\gamma}^{-1}(\widehat{x})
E_s (\widehat{x})$ qui vérifie $\|
\Phi_{0}(Y_1) - \Phi_{0}(Y_2) \| \leq \gamma_0 \| Y_1 - Y_2 \|$. Son
  image par $\tau_x \circ C_{\gamma}(\widehat{x})$ est
  la variété stable approchée que l'on voulait construire au point $x$. 

On va maintenant minorer le volume de cette variété.

\subsection{{\bf Minoration du volume des
    variétés stables}}

Dans un premier temps, on suppose que $n_1 >0$.

On repart du graphe $(\Phi_n(Y),Y)$ et cette fois-ci on va tirer en
arrière des tranches de celui-ci. Plus précisément considérons:
$$\{ 0 \}^{s-l-1} \times \{ a_{s-l} \} \times \dots \times \{
a_{k-n_1} \} \times B_1(0,e^{-3p \gamma n })$$
avec $(a_{s-l} , \dots , a_{k-n_1}) \in B_2(0 , e^{-4p \gamma n})$.

Cet élément est un graphe $(\Psi_n(Z),Z)$ au-dessus d'une partie de $C_{\gamma}^{-1}(\widehat{f}^n(\widehat{x})) E_s^1 (
\widehat{f}^n(\widehat{x}))$.

Pour les mêmes raisons que dans le paragraphe précédent, l'image de ce
graphe par $g_{\widehat{f}^n(  \widehat{x})}^{-1}$ est un graphe $(\Psi_{n-1}(Z),Z)$ au-dessus
  d'une partie de $C_{\gamma}^{-1}(\widehat{f}^{n-1}(\widehat{x}))
E_s^1 ( \widehat{f}^{n-1}(\widehat{x}))$. Il vérifie de plus $\|
  \Psi_{n-1}(Z_1) - \Psi_{n-1}(Z_2) \| \leq \gamma_0 \| Z_1 - Z_2 \|$.

Par ailleurs, ce graphe vérifie aussi le (pour $\gamma_0 << \gamma$):

\begin{lemme}

La projection du graphe $(\Psi_{n-1}(Z),Z)$ sur $C_{\gamma}^{-1}(\widehat{f}^{n-1}(\widehat{x}))
E_s^1 ( \widehat{f}^{n-1}(\widehat{x}))$ contient la boule $B_1(0, e^{-
  3 p \gamma n})$ et $\| \Psi_{n-1}(0) \| \leq e^{- 4 p\gamma n +2 \gamma}$.

\end{lemme}

\begin{proof}

On applique la deuxième partie du théorème de la transformée de graphe du paragraphe
\ref{graphe}.

Ici $\alpha=e^{-3 p \gamma n}$, $\beta= e^{-4 p \gamma
  n}$ et $\|B^{-1} \|^{-1} \geq e^{\gamma}$ car on a sélectionné les
exposants strictement négatifs.

Maintenant, on peut prendre $r=e^{-2 \gamma n p}$ (car $(\Phi_n(Y),Y)$ vit dans $B(0, e^{-2 \gamma n p})$) et en considérant les estimées obtenues à
la fin de la démonstration de la proposition \ref{inverse}, on en
déduit (pour $n$ grand):
$$ \| D^2 g_{\widehat{f}^n(  \widehat{x})}^{-1} \|_{B(0,r)} \leq 
e^{4 \gamma n} e^{\gamma n p}.$$

Comme dans le paragraphe précédent, le $\delta$ du théorème de la
transformée de graphe est très petit. On en
déduit donc bien que la projection du graphe $(\Psi_{n-1}(Z),Z)$ sur $C_{\gamma}^{-1}(\widehat{f}^{n-1}(\widehat{x}))
E_s^1 ( \widehat{f}^{n-1}(\widehat{x}))$ contient la boule $B_1(0, e^{-
  3 p \gamma n})$.

Pour la majoration de $\| \Psi_{n-1}(0) \|$ il y a deux cas. Soit
$n_1=k-s+l+1$ (i.e. $\chi_{s-l}, \dots , \chi_k$ sont tous strictement
négatifs) et alors $\| \Psi_{n-1}(0) \| = \| \Phi_{n-1}(0) \|=0$. Soit $n_1 < k-s+l+1$ et la norme $\| A \|$ du théorème du paragraphe \ref{graphe} est
majorée par $e^{\gamma}$. L'estimée donnée dans ce théorème donne donc
$\| \Psi_{n-1}(0) \| \leq e^{- 4 p\gamma n +2 \gamma}$ si $\gamma_0$
est très petit devant $\gamma$.

\end{proof}

Maintenant, on ne garde que la partie de ce graphe qui se trouve
au-dessus de $\{ 0 \}^{k-n_1} \times B_1(0, e^{-3 p \gamma
  n})$. Remarquons que ce cut-off est finalement le même que
celui du paragraphe précédent (où on gardait la partie au-dessus de $\{ 0 \}^{s-l-1} \times B_2(0, e^{- 3
  \gamma p n}) \times B_1(0,e^{- 3 p \gamma n}) $). En effet, si on prend un point $Z$ dans
$B_1(0, e^{-3 p \gamma  n})$, on a $\| \Psi_{n-1}(Z) - \Psi_{n-1}(0) \| \leq
\gamma_0 \|Z \| \leq \gamma_0 e^{-3p \gamma n}$. On en déduit que $\|
\Psi_{n-1}(Z) \| \leq 2 \gamma_0 e^{-3p \gamma n}$ et donc que la projection de $(\Psi_{n-1}(Z),Z)$ sur $C_{\gamma}^{-1}(\widehat{f}^{n-1}(\widehat{x}))
E_s ( \widehat{f}^{n-1}(\widehat{x}))$ est incluse dans $B_2(0, e^{- 3  \gamma p n}) \times B_1(0,e^{- 3 p \gamma n}) $.

On recommence maintenant tout ce que l'on vient de faire en poussant en avant par
$g_{\widehat{f}^{n-1}(\widehat{x})}^{-1}$ et ainsi de suite. A la fin,
on obtient un graphe $(\Psi_0(Z),Z)$ au-dessus d'une partie de
$C_{\gamma}^{-1}(\widehat{x}) E_s^1(\widehat{x})$ qui contient $ B_1(0, e^{- 3
  \gamma p n})$ et avec $\| \Psi_0(0) \| \leq e^{-4 \gamma p n +
  2 \gamma n}$. De plus ce graphe est assez plat car il
vérifie $\| \Psi_{0}(Z_1) - \Psi_{0}(Z_2) \| \leq \gamma_0 \|Z_1 - Z_2 \|$.

En faisant varier $(a_{s-l} , \dots , a_{k-n_1}) \in B_2(0 , e^{-4p
  \gamma n})$, on a donc feuilleté la variété stable
  approchée. Grâce à cette propriété, nous allons pouvoir minorer le
  volume $k-s+l+1$-dimensionnel de cette variété.

On se place dans $C_{\gamma}^{-1}(\widehat{x}) E_u(\widehat{x}) \oplus
C_{\gamma}^{-1}(\widehat{x}) E_s^2(\widehat{x}) \oplus
C_{\gamma}^{-1}(\widehat{x}) E_s^1(\widehat{x})$ et on considère un
plan complexe de dimension $k-n_1$ de la forme
$x_{k-n_1+1}=b_{k-n_1+1}, \dots , x_k=b_k$ avec $(b_{k-n_1+1}, \dots ,
b_k) \in B_1(0, e^{-3 \gamma p n})$. L'intersection $I_0$ de ce plan avec le
graphe $(\Phi_0(Y), Y)$ de la variété stable est de dimension
$k-n_1-s+l+1$. Nous allons minorer le volume
$k-n_1-s+l+1$-dimensionnel de cette intersection par environ $e^{-2
  \chi_{s-l}^{+} n - \dots - 2 \chi_k^{+} n }$. Cela
impliquera que le volume $k-s+l+1$-dimensionnel du graphe $(\Phi_0(Y), Y)$ (on note
$W_0$ cette variété) sera supérieur à
$$ e^{-2  \chi_{s-l}^{+} n - \dots - 2 \chi_k^{+} n - 6 \gamma p    n
  n_1}.$$
En effet, par la formule de la coaire (voir \cite{Fe} p. 258),  ce
  volume est supérieur à:
$$ \int_{B_1(0 , e^{- 3  \gamma p n})} \int_{\pi_2^{-1}(Z) \cap W_0} d \Hcal^{2(k-s+l+1-n_1)} d
    \Hcal^{2n_1}(Z),$$
(où $\pi_2$ est la projection orthogonale sur
    $C_{\gamma}^{-1}(\widehat{x}) E_s^1(\widehat{x})$). Et cette
    dernière quantité est plus grande que $e^{-2  \chi_{s-l}^{+} n - \dots - 2 \chi_k^{+} n - 6 \gamma p
    n  n_1}$. 

Il reste donc à montrer que le volume $i_0=k-n_1-s+l+1$-dimensionnel
de $I_0$ est supérieur à environ $e^{-2  \chi_{s-l}^{+} n - \dots - 2
  \chi_k^{+} n }$.

Avant cela faisons une remarque. On a supposé jusqu'ici que $n_1
>0$. Lorsque $n_1=0$, on ne fait pas le tranchage du début de ce
paragraphe et on passe directement à la minoration de $I_0$ qui est
égal à $W_0$ que l'on va faire maintenant.

Dans un premier temps, on va évaluer le volume de
$g_{\widehat{x}}(I_0)$ en fonction de celui de $I_0$. On a:
$$\int_{I_0}  \| \Lambda^{ i_0} D g_{\widehat{x}}(Z) \|^2 d
\Hcal^{2i_0}=\mbox{volume}(g_{\widehat{x}}(I_0))$$
toujours par la formule de la coaire (ici on considère
$g_{\widehat{x}}: I_0 \rightarrow  g_{\widehat{x}}(I_0)$ comme
fonction sur $I_0$ et $D g_{\widehat{x}}(Z)$ désigne toujours la
différentielle complexe). Il s'agit de majorer $\| \Lambda^{ i_0} D g_{\widehat{x}}(Z) \|$ avec
$Z \in I_0$. Tout d'abord cette quantité est inférieure à $\|
\Lambda^{ i_0} D g_{\widehat{x}}(Z) \|$ où cette fois-ci
$g_{\widehat{x}}$ est considérée comme fonction sur $W_0$. Maintenant,
$$\| \Lambda^{ i_0} D g_{\widehat{x}}(Z) \| = \| D g_{\widehat{x}}(Z)
v_1 \wedge \dots \wedge  D g_{\widehat{x}}(Z) v_{i_0} \|$$
pour certains $v_1 , \dots , v_{i_0}$ tangents à $W_0$ (voir \cite{Ar} p. 119-120 pour les propriétés des
produits extérieurs). Soient $u_1, \dots ,
u_{i_0}$ la projection sur $C_{\gamma}^{-1}(\widehat{x})
E_s(\widehat{x})$ de $v_1, \dots , v_{i_0}$. Cela signifie que pour
$i=1, \dots , i_0$ on a $v_i= (D \Phi_0(P) \alpha_i,
\alpha_i)$ et $u_i=(0 , \alpha_i)$ où $P$ est la projection
orthogonale de $Z$ sur $C_{\gamma}^{-1}(\widehat{x})
E_s(\widehat{x})$. D'après le contrôle sur le graphe $(\Phi_0(Y),Y)$, on a $\| D \Phi_0(P) \alpha_i \|\leq \gamma_0 \| \alpha_i \|$. Cela
implique que $\| v_1 \wedge \dots \wedge v_{i_0} -  u_1 \wedge \dots \wedge
u_{i_0}\| \leq \epsilon(\gamma_0)$ avec $\epsilon(\gamma_0)$ aussi
petit que l'on veut pourvu que $\gamma_0$ le soit. En effet, si on
note $G$ l'application linéaire 
$$G(X,Y)=(D \Phi_0(P) Y,Y)=
\left(
\begin{array}{cc}
0 & D \Phi_0(P)\\
0 & I
\end{array}
\right) \left(
\begin{array}{c}
X\\
Y
\end{array}
\right), $$
on a 
$$\| v_1 \wedge \dots \wedge v_{i_0} - u_1 \wedge \dots \wedge
u_{i_0}\| \leq \| (\Lambda^{i_0} G- \Lambda^{i_0} \Ical)(  u_1 \wedge \dots \wedge
u_{i_0})\|$$
où 
$$\Ical(X,Y)=\left(
\begin{array}{cc}
0 & 0\\
0 & I
\end{array}
\right) \left(
\begin{array}{c}
X\\
Y
\end{array}
\right).$$
Enfin $\| \Lambda^{i_0} G- \Lambda^{i_0} \Ical \|$ est aussi petit que l'on veut pourvu que
$\gamma_0$ le soit (par le théorème des accroissements finis).

Maintenant,
$$\| \Lambda^{ i_0} D g_{\widehat{x}}(Z) \| = \| D g_{\widehat{x}}(Z)
v_1 \wedge \dots \wedge  D g_{\widehat{x}}(Z) v_{i_0} \|$$
vérifie:
$$\| \Lambda^{ i_0} D g_{\widehat{x}}(Z) \| \leq \|\Lambda^{ i_0} D
g_{\widehat{x}}(0)(u_1 \wedge \dots \wedge u_{i_0}) \| +A+B$$
avec
$$A= \|(\Lambda^{ i_0} Dg_{\widehat{x}}(Z) - \Lambda^{ i_0} D
g_{\widehat{x}}(0))(v_1 \wedge \dots \wedge v_{i_0}) \|$$
$$B=\|\Lambda^{ i_0} D
g_{\widehat{x}}(0)(v_1 \wedge \dots \wedge v_{i_0}-u_1 \wedge \dots
\wedge u_{i_0}) \|.$$

$A$ est inférieur à  $\|\Lambda^{ i_0} Dg_{\widehat{x}}(Z) - \Lambda^{ i_0} D
g_{\widehat{x}}(0)\|$ qui est aussi petit que l'on veut pourvu que
l'on prenne $n$ grand. En effet d'une part comme à
chaque étape on fait un cut-off, $Z$ vit ici dans la boule centrée en
$(0,0)$ et de rayon $e^{- 2 \gamma p n}$, d'autre part on a un
contrôle de la différentielle seconde qui est donnée par la
proposition \ref{direct}.

Comme $B$ est inférieur à $\epsilon(\gamma_0)$, on obtient pour $n$
grand (en particulier devant des constantes qui dépendent de
$\gamma_0$): 
$$\| \Lambda^{ i_0} D g_{\widehat{x}}(Z) \| \leq  \|\Lambda^{ i_0} D
g_{\widehat{x}}(0)_{| C_{\gamma}^{-1} E_s(\widehat{x})} \| + \epsilon(\gamma_0).$$
Mais 
$$ \|\Lambda^{ i_0} D
g_{\widehat{x}}(0)_{| C_{\gamma}^{-1} E_s(\widehat{x})} \| \leq e^{ \chi_{s-l} + \dots +  \chi_{k-n_1} 
  + \gamma k}$$
ce qui implique que
$$\| \Lambda^{ i_0} D g_{\widehat{x}}(Z) \| \leq e^{ \chi_{s-l}^{+} + \dots +  \chi_{k}^{+} 
  + \gamma k}(1 + \epsilon(\gamma_0)).$$

Le volume de $g_{\widehat{x}}(I_0)$ est donc majoré par
$\mbox{volume}(I_0) \times e^{2 \chi_{s-l}^{+} + \dots + 2 \chi_{k}^{+} 
  +2 \gamma k}(1 + \epsilon(\gamma_0))^2$.

Maintenant, on prend l'image de $g_{\widehat{x}}(I_0)$ par
$g_{\widehat{f}(\widehat{x})}$ et ainsi de suite. En faisant les mêmes
calculs que précédemment, on obtient alors une majoration du volume
$i_0$-dimensionnel de $g_{\widehat{x}}^n(I_0)$ par $\mbox{volume}(I_0)
\times e^{2 \chi_{s-l}^{+}n + \dots + 2 \chi_{k}^{+} n
  + 2 \gamma kn}(1 + \epsilon(\gamma_0))^{2n}$.

Mais $g_{\widehat{x}}^n(I_0)$ rencontrent tous les $$\{ 0 \}^{s-l-1} \times \{ a_{s-l} \} \times \dots \times \{
a_{k-n_1} \} \times B_1(0,e^{-3p \gamma n })$$
avec $(a_{s-l} , \dots , a_{k-n_1}) \in B_2(0 , e^{-4p \gamma
  n})$. Le volume $i_0$-dimensionnel de
$g_{\widehat{x}}^n(I_0)$ est donc supérieur à $e^{-8p i_0 \gamma
  n}$. Autrement dit, on a minoré le volume de $I_0$ par:
$$e^{-2 \chi_{s-l}^{+}n - \dots - 2 \chi_{k}^{+} n
  - 2 \gamma kn - 8p \gamma i_0 n}(1 + \epsilon(\gamma_0))^{-2n}$$
qui est supérieur à 
$$e^{-2 \chi_{s-l}^{+}n - \dots - 2 \chi_{k}^{+} n
   - 10p \gamma k n} $$
(car $\gamma_0$ est très petit devant $\gamma$).

C'est la minoration que l'on cherchait car on obtient une minoration
du volume de $W_0$ par
$$e^{-2 \chi_{s-l}^{+}n - \dots - 2 \chi_{k}^{+} n
   - 16p \gamma k n}.$$

On va passer à la majoration du volume de toutes ces variétés stables
à l'aide des degrés dynamiques de $f$.

\subsection{{\bf Majoration du volume}}

Nous avons construit des variétés stables pour chaque $x_i$ ($i=1,
\dots , N$) au-dessus de $C_{\gamma}^{-1}(\widehat{x_i}) E_s( \widehat{x_i})$. Considérons
maintenant l'image de ces variétés par les
$C_{\gamma}(\widehat{x_i})$. Chaque image est un graphe au-dessus de
$E_s(\widehat{x_i})$ (pour le repère $E_u(\widehat{x_i}) \oplus
E_s(\widehat{x_i})$). De plus, si $(\Phi(Y),Y)$ est l'un d'eux, on a
$\| \Phi(Y_1) - \Phi(Y_2) \| \leq \frac{\gamma_0}{\alpha_0^2} \| Y_1 -
Y_2 \|$ qui est
aussi petit que l'on veut pourvu que $\gamma_0$ soit petit par rapport
à $\alpha_0$ (ce que l'on a supposé). Quitte à remplacer $N$ par $N/K$
où $K$ est une constante qui ne dépend que de $X$, on peut supposer
que tous ces graphes vivent dans une carte fixée $\psi: U \rightarrow
X$ et que les $x_i$ sont à distance au moins $\epsilon_0$ du bord de
$U$ (cela signifie que $\tau_{x_i}$ est égal à $\psi$ modulo une
translation). Toujours quitte à remplacer $N$ par $N/K$, on peut
supposer que les graphes précédents sont des graphes au-dessus d'un
plan complexe $P$ de dimension $k-s+l+1$ et que la projection de
chaque graphe sur $P$ est de volume supérieur à 
$$e^{-2 \chi_{s-l}^{+}n - \dots - 2 \chi_{k}^{+} n
   - 16p \gamma k n }$$
(éventuellement redivisé par une constante). Dans ce qui précède le
plan $P$ peut être bougé un petit peu.

Le volume $k-s+l+1$-dimensionnel de la projection de tous les graphes
$W(x_i)$ sur $P$ est donc supérieur à:

$$e^{h(\mu)n -2 \chi_{s-l}^{+}n - \dots - 2 \chi_{k}^{+} n
   - 20p \gamma k n }.$$

Nous allons maintenant majorer ce volume à l'aide des degrés
dynamiques.

Notons, $\pi_3$ la projection orthogonale sur $P$. $\pi_3(U)$ vit dans
un compact $K$ de $P$ et pour $a \in K$, $F_a$ désignera la fibre
$\pi_3^{-1}(a)$. Elle est de dimension $s-l-1$. Dans la suite $da$
sera la mesure de Lebesgue sur un voisinage de $K$ dans $P$.

Si $\Wcal_s = \cup_{i=1}^{N} W(x_i)$ et $n(a)$ désigne le nombre
d'intersection de $F_a$ avec $\Wcal_s$, on a:

$$\int n(a) da = \mbox{volume de la projection de } \Wcal_s \mbox{ sur
}P$$
est supérieur à
$$e^{h(\mu)n -2 \chi_{s-l}^{+}n - \dots - 2 \chi_{k}^{+} n
   - 20p \gamma k n }.$$

Cependant, pour $a$ fixé, les images par $\psi$ des points
d'intersection entre $F_a$ et $\Wcal_s$ sont $(n,
\delta/2)$-séparés. En effet considérons $y_1 \in \psi(F_a \cap
W(x_i))$ et $y_2 \in \psi(F_a \cap
W(x_j))$. Par définition des $x_i$, on a $d_n(x_i,x_j) \geq
\delta$. Cela signifie qu'il existe $l$ compris entre $0$ et $n-1$
avec $d(f^l(x_i),f^l(x_j)) \geq \delta$. Mais le diamètre de
$f^l(\psi(W(x_i))$ et de $f^l(\psi(W(x_j))$ est inférieur à $\delta/4$
(car on a fait des cut-off) ce qui implique que $d(f^l(y_1), f^l(y_2))
\geq \delta/2$. Les points $y_1$ et $y_2$ sont $(n, \delta/2)$-séparés.

Si $\Omega_f= X \setminus \cup_{n \in \Zz} f^n(I_f)$, on notera
$\Gamma_n(a)$ l'adhérence de $\{ (z,f(z), \dots ,
f^{n-1}(z) ) \mbox{ , } z \in \psi(F_a \cap U) \cap \Omega_f \}$ dans $X^n$. C'est le
multigraphe de $\psi(F_a \cap U)$. On munit $X^n$ de la forme de Kähler
$\omega_n= \sum_{i=1}^{n} \Pi_i^{*} \omega$ où les $\Pi_i$ sont les
projections de $X^n$ sur ses facteurs. Maintenant, on a:

\begin{lemme}

$$\int \mbox{volume}(\Gamma_n(a))da \geq c(\delta) \int n(a) da.$$

\end{lemme}

\begin{proof}

La démonstration est la même que dans \cite{BD} (paragraphe 5) et \cite{Gr2}. Elle repose sur le théorème de Lelong
(voir \cite{Le}). Nous la donnons par confort pour le lecteur.

On fixe $a$. Les $n$-orbites des points d'intersection entre
$\psi(F_a \cap U)$ et $\psi(\Wcal_s)$ induisent un ensemble $F$ de
$\Gamma_n(a)$ qui est $\delta/2$-séparé pour la métrique produit de
$X^n$. Cela signifie que les $n(a)$ boules $B(y, \delta/4)$ avec $y
\in F$ sont disjointes. Par le théorème de Lelong, le volume de $\Gamma_n(a)
\cap B(y, \delta/4)$ est minoré par une constante $c(\delta)$. On en
déduit donc que le volume de $\Gamma_n(a)$ est plus grand que
$c(\delta) n(a)$. Cela démontre le lemme.
  
\end{proof}

Maintenant,

$$\int \mbox{volume}(\Gamma_n(a))da = \int \int_{\Gamma_n(a)}
w_n^{s-l-1} da$$
est égal à
$$ \int \left( \sum_{0 \leq n_1,
  \dots , n_{s-l-1 } \leq n-1} \int_{\psi(F_a \cap U) \cap \Omega_f}
  (f^{n_1})^{*} \omega \wedge \dots \wedge (f^{n_{s-l-1}})^{*} \omega
  \right) da$$
par définition de $\omega_n$ et le fait que l'on peut prendre $P$
  générique. Soit $\Omega= \int [\psi(F_a \cap U)] da$ avec $[\psi(F_a \cap U)]$ le
courant d'intégration sur $\psi(F_a \cap U)$. $\Omega$ est une forme de bidimension $(s-l-1,s-l-1)$.

De l'égalité précédente, on déduit que:
$$\int \mbox{volume}(\Gamma_n(a))da = \sum_{0 \leq n_1,
  \dots , n_{s-l-1 } \leq n-1} \int_{\Omega_f} \Omega \wedge
  (f^{n_1})^{*} \omega \wedge \dots \wedge (f^{n_{s-l-1}})^{*}
  \omega,$$
qui est inférieur à 
$$ C_0 \sum_{0 \leq n_1,
  \dots , n_{s-l-1 } \leq n-1} \int_{\Omega_f} \omega^{k-s+l+1} \wedge
  (f^{n_1})^{*} \omega \wedge \dots \wedge (f^{n_{s-l-1}})^{*}
  \omega,$$
où $C_0$ est une constante telle que $\Omega \leq C_0 \omega^{k-s+l+1}$. En
  utilisant le lemme \ref{lemme1} avec $q=s-l-1$, on obtient:
$$\int \mbox{volume}(\Gamma_n(a))da \leq c_{\epsilon} n ^{s-l-1} (\max_{0 \leq j \leq s-l-1}
  d_j + \epsilon)^n.$$
Finalement, en combinant les inégalités obtenues, on a:
$$c_{\epsilon} n ^{s-l-1} (\max_{0 \leq j \leq s-l-1}
  d_j + \epsilon)^n \geq e^{h(\mu)n -2 \chi_{s-l}^{+}n - \dots - 2 \chi_{k}^{+} n
   - 20p \gamma k n },$$
qui implique la première inégalité du théorème.

\section{{\bf Démonstration de la deuxième inégalité du théorème}}

Comme la démonstration est à peu près la même que pour la première
formule, on ne fera que l'esquisser.

Rappelons que l'on a $\chi_1 \geq \dots \geq \chi_s = \dots =
\chi_{s+l'} > \chi_{s+l'+1}  \geq \dots \geq \chi_k$.
Pour $\widehat{x} \in \widehat{Y}_{\alpha_0}$, on notera $E_1(\widehat{x}) , \dots , E_{m}(\widehat{x})$ les $E_i(\widehat{x})$
correspondant aux exposants $\chi_1 , \dots , \chi_{s+l'}$ et
$E_{m+1}(\widehat{x}) , \dots , E_{q}(\widehat{x})$ les
$E_i(\widehat{x})$ de $\chi_{s+l'+1} , \dots , \chi_k$. Soit:
$$E_u (\widehat{x}) = \oplus_{i=1}^{m} E_i(\widehat{x}) \mbox{ et }
E_s (\widehat{x}) = \oplus_{i=m+1}^{q} E_i(\widehat{x}).$$
Par ailleurs, $E_u (\widehat{x})$ sera dans la suite coupé en deux parties. Soit $n_1$ le nombre d'exposants parmi
$\chi_{1}, \dots , \chi_{s+l'}$ qui sont strictement positifs (bien sûr
$n_1$ peut être égal à $0$). Alors, nous noterons $E_u^1(
\widehat{x})$ la somme directe des $E_i(\widehat{x})$ ($i=1, \dots,
m$) correspondant aux $\chi_i$ strictement positifs et
$E_u^2(\widehat{x})$ la somme directe des autres $E_i(\widehat{x})$ de
$E_u (\widehat{x})$. La dimension de $E_u^1(\widehat{x})$ est donc
$n_1$ et celle de $E_u^2(\widehat{x})$ est $s+l'-n_1$.

On reprend les $N$ points $x_i$ du paragraphe précédent. En chaque
$f^n(x_i)$, on peut construire des variétés instables approchées de
dimension $s+l'$ par le procédé suivant. Soit $x$ un des $x_i$ (on a $x=\pi(\widehat{x})$ avec $\widehat{x} \in \widehat{Y}_{\alpha_0}$). On se place dans $C_{\gamma}^{-1}(\widehat{x})  E_u (
\widehat{x}) \oplus C_{\gamma}^{-1}(\widehat{x}) E_s ( \widehat{x})$ et on part de 
$$ B_1(0, e^{-3p \gamma n}) \times
B_2(0, e^{-4p \gamma n}) \times \{ 0 \}^{k-s-l'},$$
où $B_1(0, e^{-3p \gamma n})$ est la boule de $\Cc^{n_1}$ de centre $0$ et de
rayon $e^{-3p \gamma n}$ et $B_2(0, e^{-4p \gamma n})$ celle de
$\Cc^{s+l'-n_1}$ de centre $0$ et de rayon $e^{-4p \gamma n}$. Cet ensemble est un graphe $(X, \Phi_0(X))$
au-dessus d'une partie de $C_{\gamma}^{-1}(\widehat{x}) E_u (
\widehat{x})$ (avec $\Phi_0(X)=0$). Toujours grâce à la transformée de
graphe et le procédé de cut-off appliqués aux $g_{\widehat{f}^i(\widehat{x})}$, on obtient un graphe $(X, \Phi_n(X))$ au-dessus
d'une partie de $C_{\gamma}^{-1}(\widehat{f}^n(\widehat{x})) E_u
(\widehat{f}^n(\widehat{x}))$. Par les mêmes arguments qu'au
paragraphe précédent, le volume $s+l'$-dimensionnel de ce
graphe est supérieur à  
$$e^{2 \chi_{1}^{-}n + \dots + 2 \chi_{s+l'}^{-} n
   - 16p \gamma k n }.$$
Pour chaque $x_i$, on considère l'image du graphe construit au-dessus
de $C_{\gamma}^{-1}(\widehat{f}^n(\widehat{x_i})) E_u
(\widehat{f}^n(\widehat{x_i}))$ par
$C_{\gamma}(\widehat{f}^n(\widehat{x_i}))$. On notera $W(x_i)$ cette
image et $\Wcal_u = \cup_{i=1}^{N} W(x_i)$. Comme dans le paragraphe
précédent, quitte à changer $N$ en $N/K$, on peut supposer que les $N$
variétés $W(x_i)$ vivent dans une carte $\psi:U \rightarrow X$ fixée, que les $W(x_i)$
sont des graphes au-dessus d'un plan $P$ de dimension $s+l'$ et que le
volume $s+l'$-dimensionnel de la projection par $\pi_4$ des $W(x_i)$ sur $P$ est
supérieur à 
$$e^{2 \chi_{1}^{-}n + \dots + 2 \chi_{s+l'}^{-} n
   - 16p \gamma k n }.$$

Maintenant, $\pi_4(U)$ vit dans un compact $K$ de $P$ et pour $a \in
K$ on notera $F_a$ la fibre $\pi_4^{-1}(a)$. Si $n(a)$ désigne le nombre d'intersection entre
$F_a$ et $\Wcal_u$ et $da$ la mesure de Lebesgue sur un voisinage de
$K$ dans $P$, on a:
$$\int n(a) da \geq e^{h(\mu)n + 2 \chi_{1}^{-}n + \dots + 2 \chi_{s+l'}^{-} n
   - 20 \gamma k n}.$$
Les points d'intersection entre $F_a$ et $\Wcal_u$ induisent un
ensemble $(n, \delta/2)$-séparé dans $f^{-n}(\psi(F_a \cap U))$. En effet les
$g_{\widehat{f}^n(\widehat{x_i})}^{-j}(C_{\gamma}^{-1}(\widehat{f}^n(\widehat{x_i}))W(x_i))$
sont de diamètre très petit (pour $j=0, \dots , n$) et les $x_i$ sont $(n, \delta)$-séparés. 

Si on note $\Gamma_n(a)$ le multigraphe de $f^{-n}(\psi(F_a \cap U))$, on a
alors (toujours par le théorème de Lelong):
$$\int \mbox{volume}(\Gamma_n(a))da \geq c(\delta) \int n(a) da.$$
Pour finir il reste à majorer $\int \mbox{volume}(\Gamma_n(a))da $.

Par un raisonnement équivalent à celui du paragraphe précédent, cette
intégrale est inférieure à:

$$ C_1 \sum_{0 \leq n_1,
  \dots , n_{k-s-l' } \leq n-1} \int_{\Omega_f} (f^n)^{*} \omega^{s+l'} \wedge
  (f^{n_1})^{*} \omega \wedge \dots \wedge (f^{n_{k-s-l'}})^{*}
  \omega,$$
où $C_1$ est une constante qui ne dépend que de $X$.
En utilisant le lemme \ref{lemme2}, on obtient:
$$c_{\epsilon} n^{k-s-l'} ( \max_{s+l' \leq j \leq k}d_j + \epsilon)^n
\geq e^{h(\mu)n + 2 \chi_{1}^{-}n + \dots + 2 \chi_{s+l'}^{-} n
   - 20 \gamma k n }.$$
Cela démontre la deuxième inégalité.

\section{\bf{Le cas des difféomorphismes de classe $C^{1+ \alpha}$}}

Dans ce paragraphe, nous suivons la demande du referee en
donnant une version de notre théorème pour les difféomorphismes de
classe $C^{1 + \alpha}$ dans les variétés Riemanniennes compactes. Nous
aboutirons ainsi à une inégalité plus faible que celle de J. Buzzi (voir \cite{Bu}). Commençons par préciser le cadre de ce
paragraphe.

Soit $X$ une variété Riemannienne lisse compacte de dimension $k$ et
$f$ un difféomorphisme de classe $C^{1 + \alpha}$.

J. Buzzi a introduit dans \cite{Bu} des notions d'entropie
directionnelle. Dans ce paragraphe, nous considèrerons la suivante: pour $p$ compris entre $1$ et $k$, on note 
$$ \Scal^p:= \{ \sigma: ]-1,1[^p \mapsto X \mbox{, } \sigma \mbox{ de
      classe } C^{\infty}    \}.$$

On définit le $p$-volume de $\sigma \in \Scal^p$ par la formule:

$$v_p(\sigma)= \int_{]-1,1[^p} | \Lambda^p T_x \sigma | d
\lambda(x),$$

où $d \lambda$ est la mesure de Lebesgue sur $]-1,1[^p$ et $|
    \Lambda^p T_x \sigma |$ est la norme de l'application linéaire
    $\Lambda^p T_x \sigma: \Lambda^p T_x (]-1,1[^p) \mapsto \Lambda^p
    T_{\sigma(x)} X$ induite par la métrique Riemannienne sur $X$
    (voir \cite{Ne}).

Nous désignerons par $\Scal^p(t)$ les éléments $\sigma$ de
      $\Scal^p$ pour lesquels le $p$-volume est inférieur ou égal à $t$.

L'entropie $p$-directionnelle de $f$ est alors définie par (voir \cite{Bu})

$$h_p(f):=\lim_{t \rightarrow 0} \lim_{\delta \rightarrow 0} \limsup_{n
  \rightarrow + \infty} \frac{1}{n} \sup_{\sigma \in \Scal^p(t)} \log r(\delta,n, \sigma(]-1,1[^p)).$$

Ici $r(\delta,n, \sigma(]-1,1[^p))$ est le cardinal maximal d'un
    ensemble $(n, \delta)$-séparé inclus dans $\sigma(]-1,1[^p)$. 

Alors, nous avons le:

\begin{theoreme}

Soient $\mu$ une mesure invariante, ergodique et $\chi_1 \geq \dots
\geq \chi_k$ les exposants de Lyapounov de $\mu$.

Fixons $1 \leq s \leq k$. On définit $l=l(s)$ par:

$$\chi_1 \geq \dots \geq \chi_{s-l-1} > \chi_{s-l} = \dots = \chi_s
\geq \chi_{s+1} \geq \dots \geq \chi_k,$$

où $s-l$ est égal à $1$ si $\chi_1 = \dots = \chi_s$.

Alors, on a l'inégalité suivante:

$$h(\mu) \leq h_{s-l-1}(f) + \chi_{s-l}^{+} + \dots + \chi_k^{+}$$

où $h(\mu)$ est l'entropie métrique de $\mu$ et $\chi_i^{+}= \max(\chi_i,0)$.

\end{theoreme}

La preuve de ce théorème s'obtient en faisant des modifications
mineures sur notre démonstration. Il s'agit d'utiliser l'introduction
du paragraphe 5, le paragraphe 5.1 (qui suivent des idées de
S. E. Newhouse et J. Buzzi (voir \cite{Ne} et \cite{Bu})) et enfin le
paragraphe 5.2 (qui diffère de \cite{Ne} et \cite{Bu} et où on réalise
la minoration du volume des variétés stables approchées en
les feuilletant par des sous-variétés stables). Tous les autres
paragraphes concernent les applications méromorphes et sont donc
inutiles pour la preuve de cette inégalité.

Expliquons un certain nombre des petites modifications qu'il faut
effectuer sur notre démonstration pour prouver l'inégalité ci-dessus.

Tout d'abord pour les rappels. La théorie de Pesin pour les
difféomorphismes de classe $C^{1 + \alpha}$ est bien connue (voir par
exemple \cite{KH}). Pour la preuve, on a besoin d'un analogue des
propositions \ref{direct} et \ref{inverse}. On trouvera
essentiellement la
démonstration de cet analogue dans la preuve du théorème S.3.1 de
\cite{KH}. Pour la transformée de graphe (voir le paragraphe
\ref{graphe}), il faut supposer $g$
de classe $C^{1 + a}$ (i.e. $\| D g(P) - D g(Q) \| \leq L \| P -
Q \|^{a}$) et se placer dans $\Rr^k$. La seule chose qui change dans le théorème c'est
qu'il faut remplacer $\| D^2 g \|_{B(0,r)} \beta^2$ par $L \beta^{1 +
  a}$ dans les formules.

Passons maintenant à la preuve de l'inégalité. Elle commence au début
du paragraphe 5. Il s'agit ici et dans toute la suite d'enlever les extensions naturelles et
la fonction $V$ (car on considère un difféomorphisme). Ensuite dans le
plan de la preuve, on remplace la dimension complexe par la dimension
réelle et les $2 \chi_i^{+}$ par $\chi_i^{+}$. A la fin on change la
phrase ``Enfin, dans le troisième paragraphe...'' par la suivante:
``Pour finir la démonstration, il suffit de majorer le nombre de points
$(n, \delta/2)$-séparés dans une fibre de $\pi_1$ en utilisant $h_{s-l-1}(f)$ et on
aboutit ainsi à l'inégalité annoncée.''

La démonstration continue avec le paragraphe 5.1. Ici il s'agit de tout
recopier en supposant $p$
grand par rapport à $1/ \alpha$, en utilisant l'analogue de la
proposition \ref{inverse} et en changeant $\Cc$ par $\Rr$.

Ensuite, dans le paragraphe 5.2, il faut remplacer les $2 \chi_i$ par
$\chi_i$, $\Hcal^{2(k-s+l+1-n_1)}$ par $\Hcal^{k-s+l+1-n_1}$ et $\Hcal^{2
  n_1}$ par $\Hcal^{n_1}$. Il s'agit après de considérer la
différentielle réelle et de changer
$$\int_{I_0}  \| \Lambda^{ i_0} D g_{\widehat{x}}(Z) \|^2 d
\Hcal^{2i_0}=\mbox{volume}(g_{\widehat{x}}(I_0))$$
par
$$\int_{I_0}  \| \Lambda^{ i_0} D g_{x}(Z) \| d
\Hcal^{i_0}=\mbox{volume}(g_{x }(I_0))$$
(formule de la coaire en réel). Enfin quand on parle de la
différentielle seconde de $g_{\widehat{x}}$, il faut remplacer cet
argument en utilisant le caractère Hölder de la différentielle de $f$.

\bigskip

Henry de Thélin

Université Paris-Sud (Paris 11)

Mathématique, Bât. 425

91405 Orsay

France

\end{document}